\documentclass[draft, 11pt, a4paper, sort&compress]{elsarticle}

\usepackage{amsmath}
  \usepackage{bm}
\usepackage{amssymb}
\usepackage{amsthm}
\usepackage{mathrsfs}
\usepackage{amsfonts} 
\usepackage{enumerate}
 \usepackage[protrusion=true,expansion=true]{microtype}
 \usepackage[numbers]{natbib}
 
   \usepackage[varg]{txfonts} 
\usepackage[hmarginratio=1:1,top=32mm,columnsep=20pt]{geometry}  
 \newtheorem{theorem}{Theorem}[section]

 \newtheorem{proposition}[theorem]{Proposition}
   
\theoremstyle{definition}
\newtheorem{definition}[theorem]{Definition}
 \theoremstyle{remark}
\newtheorem{remark}[theorem]{Remark}
\numberwithin{equation}{section}

 \newcommand{\eps}{\varepsilon} 
\newcommand{\abs}[1]{\left\vert#1\right\vert}

\newcommand{\norm}[1]{\Vert#1\Vert}

\newcommand{\inner}[1]{\left(#1\right)}

\makeatletter
\def\ps@pprintTitle{%
     \let\@oddhead\@empty
     \let\@evenhead\@empty
     \def\@oddfoot{\reset@font\hfil\thepage\hfil}
     \let\@evenfoot\@oddfoot}
 \makeatother

\begin{document}

\begin{frontmatter} 

\title
 {3D  hyperbolic Navier-Stokes equations in a thin strip: global well-posedness and hydrostatic limit in Gevrey space}

\author[ad1,ad2]{ Wei-Xi Li}
\ead{wei-xi.li@whu.edu.cn}
 
\address[ad1]{School of Mathematics and Statistics, Wuhan University, Wuhan
 430072, China}
 \address[ad2]{Hubei Key Laboratory of Computational Science, Wuhan University, 430072 Wuhan, China}

\author[ad3]{Tong Yang}
  \ead{matyang@cityu.edu.hk}

\address[ad3]{
	Department of Mathematics, City University of Hong Kong, Hong Kong
  }

\begin{abstract}
We consider the hyperbolic version of  three-dimensional anisotropic  Naver-Stokes equations in a thin strip and its hydrostatic limit that is a hyperbolic Prandtl type equations.  We   
 prove the global-in-time existence and uniqueness for   the two systems  and  the hydrostatic limit    when the  initial data belong to the Gevrey function space with index 2.   The proof is  based on  a  direct energy method by observing the damping effect in the systems.
  \end{abstract}

\begin{keyword}
3D hydrostatic Navier-Stokes equations,  global well-posedness, Gevrey class, hydrostatitc limit.

\MSC[2020] 35Q30, 76D03,76D10. 
 \end{keyword}
 
\end{frontmatter}

\section{Introduction and the main result}

There have been extensive studies on the well-posedness of the Prandtl type equations, while  most of them
are  concerned with the local-in-time existence and uniqueness. Compared with the local theory, the global in time property  is far from being well investigated. Here,  we   mention Xin-Zhang's work \cite{MR2020656} on global weak solutions and some recent papers \cite{MR4271962,2021arXiv210300681Y,MR4213671,MR4293727,2021arXiv211113052A,2021arXiv211112836P, MR4125518} on global  analytic or Gevrey solutions.    
 Note the above results   are obtained mainly in the two-dimensional setting so that the global well-posedness of  the three-dimensional case  remains  open. 
 
 In this paper, we aim to establish  global well-posedness theories 
 for some Prandtl type equations in the   
 three-dimensional (3D) setting. Precisely,  we will investigate the  global-in-time existence and uniqueness  of the hyperbolic version of 3D anisotropic Navier-Stokes equations and 3D   hydrostatic Navier-Stokes equations.
 The proof relies on an observation that the vertical diffusion  leads to a damping effect and
 the argument is  a direct   energy method.   
 Note that this argument does not apply to the classical Prandtl equation because of the lack of Poincar\'e inequality in the half-space.

The system of hydrostatic Navier-Stokes equations  play an   important role in the atmospheric and oceanic sciences   and it describes the large scale motion of geophysical flow    as a    limit  of Navier-Stokes equations in a thin domain where  the vertical scale is significantly smaller than the horizontal one.  By   a proper rescaling (cf. \cite{MR4149066,MR4125518,MR2563627}  for instance and references therein),   the 3D anisotropic  Navier-Stokes equations in a thin domain read
\begin{equation}\label{ans}
	\left\{
\begin{aligned}
& \big(\partial_t  + u^\eps\cdot\partial_x +v^\eps\partial_y  -\eps^2\Delta_x  -\partial_y^2\big) u^\eps + \partial_x p^\eps
=0,   &\quad (x,y)\in\mathbb R^{2}\times]0,1[, \\
& \eps^2 \big(\partial_t  +  u^\eps\cdot\partial_x +v^\eps\partial_y -\eps^2\Delta_x-\partial_y^2\big) v ^\eps + \partial_yp^\eps=0,&\quad (x,y)\in\mathbb R^{2}\times]0,1[,\\
&\partial_x\cdot u^\eps +\partial_yv^\eps =0, &\quad (x,y)\in\mathbb R^{2}\times]0,1[,
\end{aligned}\right.
\end{equation} 
where $u^\eps, v^\eps$ stand for tangential and normal components of the velocity field respectively, and
the viscosity coefficient is denoted by $\eps^2$.
 In this paper, the above system 
  is considered with the following no-slip Dirichlet boundary condition 
\begin{eqnarray*}
	u^\eps|_{y=0,1}=0,\quad  v^\eps|_{y=0,1}=0. 
\end{eqnarray*} 
By letting $\eps\rightarrow 0$, the first order approximation
 yields  the following  hydrostatic Navier-Stokes equations
\begin{equation}\label{para}
\left\{
\begin{aligned}
&\big(\partial_t+u\cdot\partial_x+v\partial_y-\partial_y^2\big) u + \partial_x p
=0,   \quad & (x,y)\in\mathbb R^{2}\times ]0,1[, \\
& \partial_yp=0,  \quad & (x,y)\in\mathbb R^{2}\times ]0,1[,\\
&\partial_x\cdot u +\partial_yv =0, \quad & (x,y)\in\mathbb R^{2}\times ]0,1[, \\
& u|_{y=0, 1}=0,  \quad  v|_{y=0,1}=0, \quad & x\in\mathbb R^{2}, \\
&u|_{t=0} =u_0^{H}, \quad & (x,y)\in\mathbb R^{2}\times ]0,1[.
\end{aligned}\right.
\end{equation}
Here,   $v$ is a scalar function  and    $u=(u_1, u_2)$ is  vector-valued,  standing for  the normal   and the tangential  velocity fields respectively.  Compared with the Navier-Stokes equations, there is no time evolution equation for the normal velocity $v$ and  the loss of tangential derivative property occurs  in the non-local term $v$. This is the main degeneracy feature of  the     Prandtl type equations.    Note that the classical Prandtl equations are considered  in the half-space:
	\begin{equation}\label{prandtl}
\left\{
\begin{aligned}
&\big(\partial_t+u\cdot\partial_x+v\partial_y-\partial_y^2\big) u + \partial_x p
=0,   \quad & (x,y)\in\mathbb R^{2}\times ]0,+\infty[, \\
&\partial_x\cdot u +\partial_yv =0, \quad & (x,y)\in\mathbb R^{2}\times ]0,+\infty[, \\
&u|_{y=0}=0, \quad v|_{y=0}=0,\quad \lim_{y\rightarrow+\infty} u=U, \quad & x\in\mathbb R^{2},\\
&u|_{t=0}=u_0^P, \quad & (x,y)\in\mathbb R^{2}\times ]0,+\infty[,
\end{aligned}\right.
\end{equation}
 where $p$ and $U$ are given by the trace of the  Euler flow on the boundary. 
 
 Other Prandtl type equations include   
  hydrostatic Euler  equations  and MHD boundary layer system. The former is an   inviscid form of \eqref{para} and the latter is  a system of Prandtl type equations on 
  velocity and magnetic fields. 
  For Prandtl type equations without structural assumption,  there are results showing that
   either analyticity or  Gevrey regularity  is sufficient for the well-posedness, cf.  \cite{MR2737223, MR1617542, MR3925144,lmy, MR4270479,MR4125518} and references therein. In particular,
   the Gevrey 2 function space is optimal for classical Prandtl equation  \cite{MR2601044}, while
   the optimal index for MHD boundary layer system remains unknown.  On the other hand,
    the analyticity is  necessary for the well-posedness of the hydrostatic  Navier-Stokes equations \eqref{para}, cf. \cite{MR2563627}. 
   
 Recently,   Aarach   \cite{2021arXiv211113052A}  and Paicu-Zhang \cite{2021arXiv211112836P}     
 studied the hyperbolic version of  2D hydrostatic Navier-Stokes equations  and established global solutions in analytic and Gevrey class 2, respectively. This shows  the hyperbolic feature 
 yields some stabilizing effect. Note that the hyperbolic version of the  hydrostatic Navier-Stokes equations can be derived  as a hydrostatic limit of the hyperbolic Navier-Stokes equations that    was proposed by   C. Cattaneo \cite{MR0032898} to avoid the non-physical property  of infinite propagation speed. 
 There have been many  results on the hyperbolic Navier-Stokes equations, cf.  \cite{MR3942552,MR2045417,MR3085225,MR3085226,MR2404054,CHR}. By performing proper change of scales as in \cite{2021arXiv211112836P}     we have the following  hyperbolic version of the anisotropic  Navier-Stokes equations \eqref{ans}:
 \begin{equation}\label{hans}
 	\left\{
\begin{aligned}
& \big(\partial_t^2+\partial_t  + u^\eps\cdot\partial_x +v^\eps\partial_y  -\eps^2\Delta_x  -\partial_y^2\big) u^\eps + \partial_x p^\eps
=0,   &\quad (x,y)\in\mathbb R^{2}\times]0,1[, \\
& \eps^2 \big(\partial_t^2+\partial_t  +  u^\eps\cdot\partial_x +v^\eps\partial_y -\eps^2\Delta_x-\partial_y^2\big) v ^\eps + \partial_yp^\eps=0,&\quad (x,y)\in\mathbb R^{2}\times]0,1[,\\
&\partial_x\cdot u^\eps +\partial_yv^\eps =0, &\quad (x,y)\in\mathbb R^{2}\times]0,1[,\\
&u^\eps|_{y=0, 1}=0,  \quad  v^\eps|_{y=0,1}=0,   \quad & x\in\mathbb R^{2},  \\
&(u^\eps, v^\eps)|_{t=0} =(u_0^\eps, v_0^\eps),\quad (\partial_t u,\partial_t v)|_{t=0}=(u_1^\eps, v_1^\eps),   \quad& (x,y)\in\mathbb R^{2}\times ]0,1[.
\end{aligned}\right.
 \end{equation} 
 Letting $\eps \rightarrow 0 $ in the above system \eqref{hans} gives    
 \begin{equation}\label{hyper}
\left\{
\begin{aligned}
&\big(\partial_t^2+\partial_t+u\cdot\partial_x+v\partial_y-\partial_y^2\big) u + \partial_x p=0,   \quad & (x,y)\in\mathbb R^{2}\times ]0,1[, \\
& \partial_yp=0,  \quad & (x,y)\in\mathbb R^{2}\times ]0,1[,\\
&\partial_x\cdot u +\partial_yv =0, \quad & (x,y)\in\mathbb R^{2}\times ]0,1[, \\
&u|_{y=0, 1}=0,  \quad  v|_{y=0,1}=0,   \quad & x\in\mathbb R^{2},  \\
&u|_{t=0} =u_0,\quad \partial_t u|_{t=0}=u_1,   \quad& (x,y)\in\mathbb R^{2}\times ]0,1[,
\end{aligned}\right.
\end{equation}
which is a hyperbolic version of the hydrostatic Navier-Stokes equations \eqref{para}.
For this system, according to  the recent work of Paicu-Zhang  \cite{2021arXiv211112836P}, we  have
  well-posedness in Gevrey function space rather than analytic space. 
   This is different from its parabolic version because analyticity is necessary for well-posedness of
   the hydrostatic Navier-Stokes equations without any structural assumption.  

Before stating the main result on the global well-posedness for the hyperbolic version
of the hydrostatic Navier-Stokes system  \eqref{hyper}, we briefly review some previous
related works on the Prandtl type equations as follows.
 
 \subsection{Classical  Prandtl equation} The well-posedness theory of  Prandtl equation  \eqref{prandtl}  has been well investigated, cf. \cite{MR3327535, MR3795028, MR3670620,
 MR3983729, MR3458159, MR3600083, MR3925144,  MR1476316, MR2601044, MR3429469, MR2849481,MR3461362, MR3284569,MR3493958, MR4055987,MR2020656, MR3710703, 2022arXiv220110139Y, MR3464051} and the references therein.   For the   2D case, under Oleinik's monotonicity condition the well-posedness theory in Sobolev space was justified in  the pioneer work by Oleinik  \cite{MR1697762}. This classical result was revisited in 
 two independent work of Alexandre-Wang-Xu-Yang \cite{MR3327535} and Masmoudi-Wong \cite{MR3385340} by using energy method.  However, the Sobolev well-poseness of 3D Prandtl 
 equation under suitable structure condition
 remains unsolved despite some attempts like  Liu-Wang-Yang \cite{MR3600083}.  If  Oleinik's monotonicity condition is violated,  the ill-posedness and
 the related  instability phenomena were  well investigated, cf. \cite{MR1476316,MR4264948,MR1761409,MR3464020,MR3566199,MR4028516,MR2601044,MR3458159} and the references therein.  Without any structural assumption, it is now well-understood that the Prandtl equation is well-posed in Gevrey class with optimal Gevrey index less or equal to $ 2$ by the instability
 analysis of   Dietert and G\'erard-Varet \cite{MR3925144} and the work on well-posedness by
 Dietert-G\'erard-Varet \cite{MR3925144} and  Li-Masmoudi-Yang \cite{lmy}.   This
 generalizes  the classical result of Sammartino-Caflisch \cite{MR1617542} in  the analytic framework. Similar well-posedness properties of hyperbolic Prandtl equations in Gevrey class were proven in \cite{2021arXiv211210450L}.      

On the other hand, in the fully nonlinear regime, Prandtl type system can be
derived from the  MHD system.  In this regime, the 
tangential magnetic field has stabilizing effort as shown in the  2D case   by Liu-Xie-Yang \cite{MR3882222}(see also  \cite{MR4342301, MR4102162} for the further generalization), where the Sobolev well-posedness theory was established without Oleinik's monotonicity condition on the velocity
field  provided the tangential magnetic field dominates. 
Without  any structural assumption, the Gevrey well-posedness was studied in \cite{MR4270479} with   Gevrey index less or equal to  $ 3/2$ that is not known to be optimal.

As for global-in-time existence of the classical Prandtl equation, there is an early work on weak solution by Xin-Zhang \cite{MR2020656}, and work
on analytic solution by Paicu-Zhang \cite{MR4271962}, cf. also  some other related work \cite{MR3461362, MR3710703, MR3464051}. Recently, in  \cite{2021arXiv210300681Y} the authors also proved the global well-posedenss property in Gevrey class 2.  On the other hand, the global analytic solution to MHD boundary layer system was obtainded  by Liu-Zhang \cite{ MR4213671} and Li-Xie \cite{MR4293727}.    
Note that all these global-in-time existence results are in 2D setting and  some suitable structural condition on the initial data is required.   Hence, the  global property of these systems in 3D setting remains unknown.

\subsection{Hydrostatic  Navier-Stokes equations  and related models}  Compared with Prandtl equation, the hydrostatic Navier-Stokes equations \eqref{para} is less being well understood. In fact, the Sobolev well-posedness of  the hydrostatic Navier-Stokes equations is still unclear. Under the convex assumption,  only the Gevrey  well-posedness has been obtained, cf. the recent work by
G\'{e}rard-Varet-Masmoudi-Vicol
 \cite{MR4149066} with Gevrey index  up to $9/8$ that seems not  to be optimal.   On the other hand, Masmoudi-Wong \cite{MR2898740} proved the convex condition is sufficient for the Sobolev well-posedness  of hydrostatic  Euler equations which is  the inviscid form of hydrostatic  Navier-Stokes equations. And M.Renardy \cite{MR2925113} obtained the classical solutions to hydrostatic MHD equations  provided the horizontal component of the magnetic field is not degenerate.    

Furthermore, the global well-posedness property of   the  hydrostatic Navier-Stokes equations \eqref{para} was investigated     by Paicu-Zhang-Zhang \cite{MR4125518} in analytic function space. 
In addition,  the global well-posedness theory of the hyperbolic version  of 2D hydrostatic Navier-Stokes equations \eqref{hyper} was established recently by Aarach   \cite{2021arXiv211113052A}  and  Paicu-Zhang \cite{2021arXiv211112836P} in analytic and Gevrey function spaces respectively.

 \subsection{Statement of the main results} In this paper, we study  the global Gevrey well-posedness  of  the  hyperbolic version \eqref{hyper} for   3D hydrostatic Navier-Stokes system. For this, we first list some
 notations to be used.
 
 \smallskip
 \noindent {\bf Notation.}  In the following, we will use  $\norm{\cdot}_{L^2}$ and $\inner{\cdot, \cdot}_{L^2}$ to denote the norm and inner product of  $L^2=L^2(\mathbb R^2\times[0,1])$   and use the notation   $\norm{\cdot}_{L_x^2}$ and $\inner{\cdot, \cdot}_{L_x^2}$  when the variable $x$ is specified. Similar notation  will be used for $L^\infty$. In addition, we use $L^p_x(L^q_y) = L^p (\mathbb R^2; L^q([0,1]))$ for the classical Sobolev space.  For a vector-valued function $A=(A_1,A_2, \ldots, A_n)$, we used the convention that $\norm{A}_{L^2}^2=\sum_{1\leq j\leq n}\norm{A_j}_{L^2}^2$
 
 \smallskip
 
   In the following discussion,  we only require the Gevrey regularity in the tangential variable $x\in\mathbb R^2$.  Precisely, the Gevrey function spaces are defined as follows.

\begin{definition} 
\label{defgev}  The space $X_{\rho}$ of  (partial) Gevrey functions  consists of all  smooth  (scalar or vector-valued) functions
 $h(t,x)$  such that the  norm  $\abs{h(t)}_{X_{\rho(t)}}<+\infty,$  where  
 \begin{eqnarray*}
\begin{aligned}
	 \abs{h}_{X_\rho  }^2= &\sum_{j=1}^2\sum_{ m= 0}^{+\infty}  L_{\rho,m}^2\Big( \norm{ \partial_t\partial_{x_j}^m     h }_{L^2}^2+ \norm{ \partial_y\partial_{x_j}^m      h }_{L^2}^2 +\norm{  \partial_{x_j}^m    h}_{L^2}^2  \Big),
	 \end{aligned}
\end{eqnarray*}
with   
 \begin{equation}\label{af}
  	L_{\rho,m}=
  	\frac{\rho^{m+1}(m+1)^7}{ (m!)^2}, \quad m\geq 0, \ \rho>0.
   \end{equation}
 In the following discussion, $\rho$ depends on time but we only write it as $\rho$ for simplicity
 of notations. 
 On the other hand,  if $h$ is independent of $t$,  then  we use the notation
\begin{eqnarray*}
\begin{aligned}
	 \abs{h}_{X_{\rho_*}  }^2= &\sum_{j=1}^2\sum_{ m=0}^{+\infty}  L_{\rho_*, m}^2\Big( \norm{ \partial_y\partial_{x_j}^m     h }_{L^2}^2 +\norm{  \partial_{x_j}^m     h }_{L^2}^2\Big)	
\end{aligned}
\end{eqnarray*}  
with $\rho_*$ being a  real number. 
\end{definition}

\begin{remark}
The norm $\abs{h}_{X_\rho}$ defined above is  equivalent to the standard Gevrey norm 
\begin{eqnarray*}
	\norm{h}_{\rho}^2=\sum_{ \alpha \geq 0}  L_{\rho,\abs\alpha}^2\Big( \norm{ \partial_t\partial_{x}^\alpha     h}_{L^2}^2+ \norm{ \partial_y\partial_{x}^\alpha    h }_{L^2}^2 +\norm{  \partial_{x}^\alpha      h }_{L^2}^2  \Big),
\end{eqnarray*}
  in the sense that
\begin{eqnarray*}
 \frac{1}{2} \norm{h}_{\rho/2}^2	\leq 	\abs{h}_{X_\rho}^2\leq   \norm{h}_{\rho}^2,
\end{eqnarray*}
 where the last inequality is trivial and the first inequality follows from  the fact that
\begin{equation}\label{equofn}
 	\forall \ \alpha\in\mathbb Z_+^2,\quad 
 	\norm{\partial_x^\alpha u}_{L^2}^2\leq  \norm{\partial_{x_1}^{\abs\alpha} u}_{L^2}^2+\norm{\partial_{x_2}^{\abs\alpha} u}_{L^2}^2.
 \end{equation}
 \end{remark}

  Note the initial data $u_0, u_1$ of \eqref{hyper} satisfy the  following  compatibility condition 
\begin{equation}\label{comp++}
 \forall\ x\in\mathbb R^2,\quad  u_0|_{y=0,1}=u_1|_{y=0,1}	=0\  \textrm{ and }\ \int_0^1 \partial_x\cdot u_0 (x, y) dy=\int_0^1 \partial_x\cdot u_1 (x, y) dy=0.
 \end{equation}
 The main results of this paper can now be stated as follows.

\begin{theorem}[Global well-posedness of system \eqref{hyper}]\label{thhyp} 
 Let $\big(X_\rho, \ \abs{\cdot}_{  X_\rho}\big)$ be given in Definition \ref{defgev}. If the initial data $u_0, u_1  \in X_{2\rho_0}$ for some $\rho_0>0$ satisfying the compatibility \eqref{comp++} and 
 \begin{eqnarray*}
 	|u_0|_{X_{2\rho_0}}  +|u_1|_{X_{2\rho_0}}\leq \eps_0
 \end{eqnarray*}
 for some small $\eps_0>0$,
  then the hyperbolic version of 3D hydrostatic Navier-Stokes equations \eqref{hyper} admit a unique global-in-time solution $u\in C\inner{[0, +\infty[, X_\rho}$ provided    
$\eps_0 
$ is sufficiently small.   Moreover 
\begin{eqnarray*}
	\forall\ t\geq 0,\quad 	\abs{u(t)}_{X_{\rho(t)}} \leq  4\eps_0 e^{-t/32} ,
\end{eqnarray*}
where and throughout the paper
\begin{equation}\label{derho}
  	\begin{aligned}
  		\rho(t)=\frac{\rho_0}{2}+\frac{\rho_0}{2} e^{-at}, \quad a=\frac{1}{96}.   
  	\end{aligned}
  \end{equation}
 If suppose additionally that 
$\partial_yu_0 \in X_{2\rho_0}$ with $|\partial_yu_0|_{X_{2\rho_0}}  \leq \eps_0$, then $\partial_t u\in C\inner{[0, +\infty[, X_\rho}$  and  $\partial_y u\in C\inner{[0, +\infty[, X_{\rho/2}}$. Moreover,  there exists a constant $C$, depending only on the Sobolev embedding constant, such that \begin{eqnarray*}
	\forall\ t\geq 0,\quad 	\abs{\partial_t u(t)}_{X_{\rho(t)}}  +\abs{\partial_y u(t)}_{X_{\rho(t)/2}} \leq  C \eps_0 e^{-t/32}.
\end{eqnarray*}
	 \end{theorem}

 \begin{theorem}[Global well-posedness of  system \eqref{hans}]\label{thmhas}
 	Suppose the initial data in \eqref{hans} satisfy that  $ (u_j^\eps,\eps v_j^\eps) \in  X_{2\rho_0}, j=0,1,$ for some $\rho_0>0$,  compatible to the boundary  conditions  in \eqref{hans}.  
  Then the anisotropic  hyperbolic Navier-Stokes equations  \eqref{hans} admit a unique global-in-time solution $(u^\eps, \eps v^\eps)\in C\inner{[0, +\infty[, X_\rho}$,  provided 
  $$
 	|(u_0^\eps, \eps v_0^\eps)|_{X_{2\rho_0}}  +|(u_1^\eps, \eps v_1^\eps)|_{X_{2\rho_0}} \leq \delta_0
$$ with $\delta_0$  sufficiently small.   Moreover,
\begin{eqnarray*}
	\forall\ t\geq 0,\quad 	\abs{\big(u^\eps (t), \eps v^\eps(t)\big)}_{X_{\rho(t)}} \leq  4\delta_0e^{-t/32},
\end{eqnarray*}
where $\rho$ is defined by \eqref{derho}.
 \end{theorem}
 
 \begin{theorem}[Hydrostatic limit]\label{thmlim}
 	Suppose all the assumptions in Theorems \ref{thhyp} and \ref{thmhas} hold, that is, the initial data of \eqref{hyper} and \eqref{hans} satisfy 
 	\begin{eqnarray*}
 	|u_0|_{X_{2\rho_0}}  +|u_1|_{X_{2\rho_0}}+|\partial_y u_0|_{X_{2\rho_0}}\leq \eps_0,
 \end{eqnarray*}
 and
 \begin{eqnarray*}
 	|(u_0^\eps, \eps v_0^\eps)|_{X_{2\rho_0}}  +|(u_1^\eps, \eps v_1^\eps)|_{X_{2\rho_0}} \leq \delta_0
 \end{eqnarray*}
 for some small constants $\eps_0$ and $\delta_0$.
  Let $u$ and $  (u^\eps,v^\eps)$ be given in Theorems  \ref{thhyp} and \ref{thmhas}  that solve   \eqref{hyper} and \eqref{hans}, respectively.   Then there exists a constant $C$, depending only on the constants $\eps_0,\delta_0,\rho_0$   and the Sobolev embedding constant but independent of $\eps$, such that
 	\begin{eqnarray*}
 		\sup_{t \geq 0}\abs{  u^\eps(t)-u(t)}_{X_{\rho(t)/2}}  \leq C\inner{\abs{u_0^\eps-u_0}_{X_{2\rho_0}}+\abs{u_1^\eps-u_1}_{X_{2\rho_0}}+  \eps}.
 	\end{eqnarray*}
 \end{theorem}
  
   \begin{remark}
  	The following analysis  implies  that the global well-posedness property  holds for  Gevrey function space with the Gevrey index less or equal to $ 2$. 
  \end{remark}

\begin{remark} The proof given in this paper is  based on  a  direct  energy method that is substantially  different from the elegant  and subtle arguments used in  \cite{2021arXiv211113052A, 2021arXiv211112836P} that  involve the  Littlewood-Paley decomposition.
 \end{remark}

 \section{Global well-posedness of hydrostatic system}\label{secderivation}
 
In this section, we will prove Theorem \ref{thhyp} that is based on  the proof of  a  priori estimates so that   the existence and uniqueness follow  by a standard argument.   In fact,  a self-contained proof  consists of two parts. 
 The first part is about the construction of approximate solutions that follows from the standard parabolic and hyperbolic theories.   And then the uniform estimate on approximate solutions can be derived following  a priori estimates.  Hence,  for brevity we only present the proof of  a priori estimates for Gevrey solutions.  
 Precisely, we will prove   the following theorem.
 
 \begin{theorem}
 	[A priori estimate]\label{thmofa}
 	If $u \in C\inner{[0,+\infty[, X_\rho} $ solves the hyperbolic version \eqref{hyper} of the  hydrostatic Navier-Stokes equations with the initial data  satisfying 
 \begin{equation}\label{asonin}
 	|u_0|_{X_{2\rho_0}} +|u_1|_{X_{2\rho_0}} \leq   \epsilon_0
 \end{equation}
 for some small small $\epsilon_0>0$, 
 then 
 \begin{equation}\label{arpe}
 		\forall\ t\geq 0,\quad 	\abs{u(t)}_{X_{\rho(t)}} \leq  4\eps_0e^{-t/32},
 \end{equation}
 where the function $\rho$ is defined by  \eqref{derho}. Moreover, suppose 
 \begin{eqnarray*}
 		|\partial_yu_0|_{X_{2\rho_0}}   \leq   \eps_0,
 \end{eqnarray*}
 and 
  $\partial_t u \in C\inner{[0,+\infty[, X_\rho} ,\partial_y u \in C\inner{[0,+\infty[, X_{\rho/2}} $. Then  there exists a constant $C$ depending only on $\rho_0$ and the Sobolev embedding constant such that
  \begin{equation}\label{secass}
	\forall\ t\geq 0,\quad 	\abs{\partial_t u(t)}_{X_{\rho(t)}} \leq  C\eps_0e^{-t/32},
\end{equation}
  and
	  \begin{equation}\label{yyu}
	\forall\ t\geq 0,\quad 	\abs{\partial_y u(t)}_{X_{\rho(t)/2}} \leq  C\eps_0e^{-t/32}.
\end{equation}
 \end{theorem}
 
 We will use a bootstrap principle to prove the above a priori estimate. To do so, we first recall an abstract version of the bootstrap principle given in \cite{MR2233925}.

 \begin{proposition}[Proposition 1.21 of  \cite{MR2233925}]\label{probootst}
Let $I$ be a time interval,  and for   each $t\in I$  we have two  statements, a ``hypothesis" $\mathbf H(t)$ and  a ``conclusion" $\mathbf C(t)$. Suppose we can verify the following four statements:  
\begin{enumerate}[(i)]
\item  If $		 \mathbf{H}( t)$ is true for some time $t\in I$ then  $\mathbf{C}( t) $  is also true for the time $t.$
		\item If   $\mathbf{C}(t)$ is true for some $t\in I$ then   $\mathbf{H}(t')$ holds for all $t'$ in a neighborhood of  $t$.
		\item If $t_1,t_2, \ldots$ is a sequence of times in $I$ which converges to another time $t\in I$ and     $\mathbf{C}(t_n)$ is true for all $ t_n$, then   $\mathbf{C}(t)$ is true.
		\item  $\mathbf{H}(t)$ is true for at least one time  $t\in I.$
	\end{enumerate}
	Then   $\mathbf{C}( t)$ is true for all $t\in I.$
\end{proposition}

  The rest  of the section is  to apply  this  bootstrap principle  to obtain the a priori estimate in Theorem \ref{thmofa}. The Gevrey class enables us to overcome the loss of tangential
  derivatives by shrinking the radius $\rho$. For this , we can either apply the abstract   Cauchy-Kowalewski Theorem (cf. \cite{lmy,MR4055987, MR1617542} for instance) or use an   auxilliary norm $\abs{\cdot}_{Y_\rho}$  as  in \cite{MR2975371}.
 Here we will make use of  the latter approach. Define  
\begin{equation}\label{ynorm}
  	\abs{h}_{Y_\rho} ^2=\sum_{j=1}^2\sum_{m=0}^{+\infty}   L_{\rho,m}^2 \Big[\frac{m+1}{\rho}  \big( \norm{ \partial_t\partial_{x_j}^{m}h}_{L^2}^2+\norm{  \partial_y\partial_{x_j}^{m}h }_{L^2}^2\big)  + \frac{(m+1)^3}{\rho^3}    \norm{ \partial_{x_j}^{m}h }_{L^2}^2 \Big]
  \end{equation}     
  for (scalar or vector-valued) functions $h$, 
 where $L_{\rho,m}$ is given by \eqref{af}.

\subsection{Proof of Theorem \ref{thmofa}: The first assertion}  In this part we will present in details the proof of    the first assertion \eqref{arpe}.   
  	To simplify the notation, we  assume without loss of generality that $\rho_0\leq 1.$ 
In the  following discussion,  we will omit the time dependence of $\rho$  in the notation, and   denote by $\rho'$ and $\rho''$ the first and the second order time derivatives respectively.   Note that
\begin{equation}\label{prpof}
\forall \ t\geq 0, \quad   \rho_0/2\leq \rho(t)\leq  \rho_0,\quad  	\rho' (t)\leq\rho'^3(t)\leq 0.
	\end{equation}    
For each $t\in [0,+\infty[$, let $\mathbf{H}(t)$	be the statement 
\begin{equation}\label{has}
 \forall \ s\in[0, t],\quad 	\abs{u(s)}_{X_{\rho(s)}} \leq 8\eps_0e^{-s/32},\end{equation}
and let  $\mathbf{C}(t)$	be the statement  
\begin{equation}\label{assc}
	\forall \ s\in[0, t],\quad 	\abs{u(s)}_{X_{\rho(s)}} \leq  4\eps_0e^{-s/32}.
\end{equation}
The statements (ii)-(iv) in Proposition \ref{probootst} follow from the continuity of $t\mapsto \abs{u(t)}_{X_{\rho(t)}}$  and the condition \eqref{asonin}.   
Then by Proposition \ref{probootst},  $\mathbf{C}(t)$ holds for all $t\in[0,+\infty[$ if we can show the
following statement: 
\begin{equation}\label{last sta}
	\begin{aligned}
		 \mathbf{H}( t) ~\textrm{ is true for some time} ~t \in [0, \infty[ \  \Longrightarrow  \mathbf{C}(t) ~\textrm{ is also true for the same time} ~t.
	\end{aligned}
\end{equation}
  We now turn to prove \eqref{last sta}. 
  In the following argument, we  assume \eqref{has} holds with some fixed time $t$.  Applying $  \partial_{x_j}^{m}, j=1,2, $ to the first equation in \eqref{hyper} gives
 \begin{equation}\label{eqmu}	
 \begin{aligned} 
	  \big(\partial_t^2+\partial_t  -\partial_y^2\big)  \partial_{x_j}^{m} u&=-\partial_x\partial_{x_j}^mp -\sum_{k=0}^{m}{m\choose k} \Big[(\partial_{x_j}^k u  \cdot\partial_{x})\partial_{x_j}^{m-k} u+(\partial_{x_j}^k v)\partial_{x_j}^{m-k}\partial_y u \Big] \\
	 &\stackrel{def}{=} H_m.
	 \end{aligned}
\end{equation}
Taking the $L^2$  inner product with 
 $ \partial_{x_j}^m u$ and $ \partial_t\partial_{x_j}^mu$, respectively, on both sides of \eqref{eqmu} and    observing $u|_{y=0,1}=0$,  we obtain
\begin{equation}\label{eqd}
	 	\frac{1}{2}\frac{d}{dt} \Big[\frac{d}{dt } \norm{ \partial_{x_j}^{m} u}_{L^2}^2+\norm{ \partial_{x_j}^{m} u}_{L^2}^2\Big]  + \norm{ \partial_y \partial_{x_j}^{m} u }_{L^2}^2=  \big( H_m, \    \partial_{x_j}^{m}u\big)_{L^2}   
	 +\norm{ \partial_t\partial_{x_j}^{m} u }_{L^2}^2,
\end{equation}
and  
\begin{equation}\label{eqt}
	 \frac{d}{dt} \inner{\norm{ \partial_t\partial_{x_j}^{m}u}_{L^2}^2+\norm{  \partial_y\partial_{x_j}^{m}u }_{L^2}^2}+ 2 \norm{ \partial_t\partial_{x_j}^{m} u}_{L^2}^2= 2\big( H_m,\  \partial_t\partial_{x_j}^{m} u\big)_{L^2}. 
\end{equation}
This yields by taking the summation of  \eqref{eqd} and \eqref{eqt}
\begin{eqnarray*}\label{fder}
	\begin{aligned}
			&\frac{1}{2}\frac{d}{dt} \Big[2\norm{ \partial_t\partial_{x_j}^{m}u}_{L^2}^2+2\norm{  \partial_y\partial_{x_j}^{m}u }_{L^2}^2+\frac{d}{dt } \norm{ \partial_{x_j}^{m} u}_{L^2}^2+\norm{ \partial_{x_j}^{m} u}_{L^2}^2\Big] +  \norm{ \partial_t\partial_{x_j}^{m} u}_{L^2}^2+ \norm{ \partial_y \partial_{x_j}^{m} u }_{L^2}^2\\
			&=  \big( H_m, \    \partial_{x_j}^{m}u+2\partial_t\partial_{x_j}^{m} u\big)_{L^2}.
			\end{aligned} 
\end{eqnarray*}
The above equality and the Poincar\'e inequality in the interval $[0,1]$ 
\begin{equation*}\label{poies}
	\frac{1}{4}\norm{\partial_{x_j}^{m} u }_{L^2}^2 \leq  \norm{ \partial_y \partial_{x_j}^{m} u }_{L^2}^2
\end{equation*}
because $u|_{y=0,1}=0$
 give
\begin{multline*}
			 \frac{1}{2}\frac{d}{dt} \Big[2\norm{ \partial_t\partial_{x_j}^{m}u}_{L^2}^2+2\norm{  \partial_y\partial_{x_j}^{m}u }_{L^2}^2+\frac{d}{dt } \norm{ \partial_{x_j}^{m} u}_{L^2}^2+\norm{ \partial_{x_j}^{m} u}_{L^2}^2\Big] \\+\frac{1}{8}\norm{  \partial_{x_j}^{m} u }_{L^2}^2+  \norm{ \partial_t\partial_{x_j}^{m} u}_{L^2}^2+\frac{1}{2} \norm{ \partial_y \partial_{x_j}^{m} u }_{L^2}^2 
			 \leq    \big( H_m, \    \partial_{x_j}^{m}u+2\partial_t\partial_{x_j}^{m} u\big)_{L^2} 	 .
\end{multline*}
Thus, by multiplying  the above inequality  by $L_{\rho,m}^2$ and observing $\frac{d}{dt}L_{\rho,m}^2= 2\rho'\frac{m+1}{\rho}L_{\rho,m}^2$,   taking summation over $m$ gives  
\begin{equation} \label{eet}
	\begin{aligned}
		& \frac{1}{2} \frac{d}{dt} \sum_{m=0}^{+\infty} L_{\rho,m}^2 \Big(2\norm{ \partial_t\partial_{x_j}^{m}u}_{L^2}^2+2\norm{  \partial_y\partial_{x_j}^{m}u }_{L^2}^2+\frac{d}{dt } \norm{ \partial_{x_j}^{m} u}_{L^2}^2+\norm{ \partial_{x_j}^{m} u}_{L^2}^2\Big) \\ &\quad+\frac18 \sum_{m=0}^{+\infty}L_{\rho,m}^2\inner{ \norm{  \partial_{x_j}^{m} u }_{L^2}^2+  \norm{ \partial_t\partial_{x_j}^{m} u}_{L^2}^2+ \norm{ \partial_y \partial_{x_j}^{m} u }_{L^2}^2} \\
			& \leq   \sum_{m=0}^{+\infty}  \rho'\frac{m+1}{\rho}L_{\rho,m}^2\Big(2\norm{ \partial_t\partial_{x_j}^{m}u}_{L^2}^2+2\norm{  \partial_y\partial_{x_j}^{m}u }_{L^2}^2+\frac{d}{dt } \norm{ \partial_{x_j}^{m} u}_{L^2}^2+\norm{ \partial_{x_j}^{m} u}_{L^2}^2\Big)   \\
			&\quad+ \sum_{m=0}^{+\infty}L_{\rho,m}^2  \big( H_m, \    \partial_{x_j}^{m}u+2\partial_t\partial_{x_j}^{m} u\big)_{L^2} 	 \\
			& \leq   \sum_{m=0}^{+\infty}  \rho'\frac{m+1}{\rho}L_{\rho,m}^2\big( \norm{ \partial_t\partial_{x_j}^{m}u}_{L^2}^2+ \norm{  \partial_y\partial_{x_j}^{m}u }_{L^2}^2  \big) + \sum_{m=0}^{+\infty}  \rho'\frac{m+1}{\rho}L_{\rho,m}^2 \norm{ \partial_t\partial_{x_j}^{m}u}_{L^2}^2   \\
			&\quad+ \sum_{m=0}^{+\infty}  \rho'\frac{m+1}{\rho}L_{\rho,m}^2 \frac{d}{dt } \norm{ \partial_{x_j}^{m} u}_{L^2}^2   + \sum_{m=0}^{+\infty}L_{\rho,m}^2\big| \big( H_m,     \partial_{x_j}^{m}u\big)_{L^2} 	\big|+2\sum_{m=0}^{+\infty}L_{\rho,m}^2\big|\big( H_m,   \partial_t\partial_{x_j}^{m} u\big)_{L^2}\big|, 
	\end{aligned}
\end{equation}
  where the last inequality holds because of 
  $\rho'\leq 0$. Moreover, for the second summation  term on the right hand side of \eqref{eet},
  by noting  $\rho'\leq 0,$ we have
\begin{eqnarray*}
\begin{aligned}
		   \rho'\frac{m+1}{\rho}L_{\rho,m}^2 \norm{ \partial_t\partial_{x_j}^{m}u}_{L^2}^2 
		 & = \rho'\frac{m+1}{\rho}  \norm{ \partial_t \big(L_{\rho,m} \partial_{x_j}^{m}u\big)}_{L^2}^2+ \rho'^3 \frac{(m+1)^3}{\rho^3}L_{\rho,m}^2 \norm{\partial_{x_j}^{m}u}_{L^2}^2\\
		 &\quad   -2 \rho'^2 \frac{(m+1)^2}{\rho^2}  \inner{ \partial_t \big(L_{\rho,m} \partial_{x_j}^{m}u\big),\ L_{\rho,m}\partial_{x_j}^{m}u}_{L^2}\\
		& \leq   \rho'^3 \frac{(m+1)^3}{\rho^3}L_{\rho,m}^2 \norm{\partial_{x_j}^{m}u}_{L^2}^2  -  \rho'^2 \frac{(m+1)^2}{\rho^2} \frac{d}{dt} \norm{ L_{\rho,m} \partial_{x_j}^{m}u }_{L^2}^2.
		 \end{aligned}
	\end{eqnarray*}
Here,  the last term can be written as   
\begin{multline*}
	 -  \rho'^2 \frac{(m+1)^2}{\rho^2} \frac{d}{dt} \norm{ L_{\rho,m} \partial_{x_j}^{m}u }_{L^2}^2 \\
	  =-  \frac{d}{dt} \Big( \rho'^2 \frac{(m+1)^2}{\rho^2}L_{\rho,m}^2\norm{    \partial_{x_j}^{m}u}_{L^2}^2\Big)+2\rho'\inner{ \rho''-\rho'^2/\rho }\frac{(m+1)^2}{\rho^2} L_{\rho,m} ^2\norm{  \partial_{x_j}^{m}u}_{L^2}^2\\
	 \leq -  \frac{d}{dt} \Big( \rho'^2 \frac{(m+1)^2}{\rho^2}L_{\rho,m}^2\norm{    \partial_{x_j}^{m}u}_{L^2}^2\Big),  
\end{multline*}
where the last inequality holds because 
\begin{equation}\label{keq}
	\rho''-\rho'^2/\rho= \frac{\rho_0 a^2 e^{-at}}{2(1+e^{-at})} \geq 0.
\end{equation}
 Combining the above inequalities yields
\begin{eqnarray*}
		 \rho'\frac{m+1}{\rho}L_{\rho,m}^2 \norm{ \partial_t\partial_{x_j}^{m}u}_{L^2}^2 
		 \leq   \rho'^3 \frac{(m+1)^3}{\rho^3}L_{\rho,m}^2 \norm{\partial_{x_j}^{m}u}_{L^2}^2 
		    -  \frac{d}{dt} \Big( \rho'^2 \frac{(m+1)^2}{\rho^2}L_{\rho,m}^2\norm{   \partial_{x_j}^{m}u}_{L^2}^2\Big) . 
\end{eqnarray*}
 For the third summation term on the right hand side of \eqref{eet}, direct computation gives
\begin{eqnarray*}
\begin{aligned}
	\rho'\frac{m+1}{\rho}L_{\rho,m}^2 \frac{d}{dt } \norm{ \partial_{x_j}^{m} u}_{L^2}^2 
	&= \frac{d}{dt }\Big(  \rho'\frac{m+1}{\rho}L_{\rho,m}^2 \norm{ \partial_{x_j}^{m} u}_{L^2}^2\Big)- \big(\rho''-\rho'^2/\rho\big)\frac{m+1}{\rho}L_{\rho,m}^2 \norm{ \partial_{x_j}^{m} u}_{L^2}^2\\  &\quad - 2  \rho'^2\frac{(m+1)^2}{\rho^2}  L_{\rho,m}^2 \norm{ \partial_{x_j}^{m} u}_{L^2}^2\\ 
	&\leq \frac{d}{dt }\Big(  \rho'\frac{m+1}{\rho}L_{\rho,m}^2 \norm{ \partial_{x_j}^{m} u}_{L^2}^2\Big) - 2  \rho'^2\frac{(m+1)^2}{\rho^2}  L_{\rho,m}^2 \norm{ \partial_{x_j}^{m} u}_{L^2}^2,
\end{aligned}	
\end{eqnarray*}
 where  \eqref{keq} is also used in the last inequality.  Now we substitute the above 
 two estimates into \eqref{eet} to have
 \begin{equation}\label{+pre}
	\begin{aligned}
		& \frac{1}{2} \frac{d}{dt} \sum_{m=0}^{+\infty} L_{\rho,m}^2 \Big(2\norm{ \partial_t\partial_{x_j}^{m}u}_{L^2}^2+2\norm{  \partial_y\partial_{x_j}^{m}u }_{L^2}^2+\frac{d}{dt } \norm{ \partial_{x_j}^{m} u}_{L^2}^2+\norm{ \partial_{x_j}^{m} u}_{L^2}^2 \Big) \\
		&\quad+  \frac{d}{dt} \sum_{m=0}^{+\infty}  L_{\rho,m}^2\norm{   \partial_{x_j}^{m}u}_{L^2}^2\Big(    \rho'^2 \frac{(m+1)^2}{\rho^2}- \rho'\frac{m+1}{\rho} \Big)\\
		 &\quad+\frac18 \sum_{m=0}^{+\infty}L_{\rho,m}^2\Big( \norm{  \partial_{x_j}^{m} u }_{L^2}^2+  \norm{ \partial_t\partial_{x_j}^{m} u}_{L^2}^2+ \norm{ \partial_y \partial_{x_j}^{m} u }_{L^2}^2 \Big) + 2\sum_{m=0}^{+\infty}  \rho'^2\frac{(m+1)^2}{\rho^2}  L_{\rho,m}^2 \norm{ \partial_{x_j}^{m} u}_{L^2}^2 \\
		&\leq \rho' \sum_{m=0}^{+\infty}  \frac{m+1}{\rho}L_{\rho,m}^2\big( \norm{ \partial_t\partial_{x_j}^{m}u}_{L^2}^2+ \norm{  \partial_y\partial_{x_j}^{m}u }_{L^2}^2  \big) + \rho'^3 \sum_{m=0}^{+\infty}  \frac{(m+1)^3}{\rho^3}L_{\rho,m}^2 \norm{\partial_{x_j}^{m}u}_{L^2}^2 
		      \\
		      &\quad+ \sum_{m=0}^{+\infty}L_{\rho,m}^2\big| \big( H_m, \    \partial_{x_j}^{m}u\big)_{L^2} 	\big|+2\sum_{m=0}^{+\infty}L_{\rho,m}^2\big|\big( H_m,\  \partial_t\partial_{x_j}^{m} u\big)_{L^2}\big|\\
			&\leq\rho'^3\sum_{m=0}^{+\infty}L_{\rho,m}^2 \Big[   \frac{m+1}{\rho}\big( \norm{ \partial_t\partial_{x_j}^{m}u}_{L^2}^2+\norm{  \partial_y\partial_{x_j}^{m}u }_{L^2}^2 \big)+     \frac{(m+1)^3}{\rho^3}  \norm{\partial_{x_j}^{m}u}_{L^2}^2\Big]\\
			&\quad+ \sum_{m=0}^{+\infty}L_{\rho,m}^2\big| \big( H_m, \    \partial_{x_j}^{m}u\big)_{L^2} 	\big|+2\sum_{m=0}^{+\infty}L_{\rho,m}^2\big|\big( H_m,\  \partial_t\partial_{x_j}^{m} u\big)_{L^2}\big|,
	\end{aligned}
\end{equation}
where \eqref{prpof} is used in the last inequality. We claim that there exists a constant $C$  depending only on  the Sobolev embedding constant such that
\begin{equation}\label{kees}
		\sum_{j=1}^2\sum_{m=0}^{+\infty}L_{\rho,m}^2\big| \big( H_m, \    \partial_{x_j}^{m}u\big)_{L^2} 	\big|+2\sum_{j=1}^2\sum_{m=0}^{+\infty}L_{\rho,m}^2\big|\big( H_m,\  \partial_t\partial_{x_j}^{m} u\big)_{L^2}\big| \leq  C   \rho^{-2}   \abs{u}_{X_{\rho}}\abs{u}_{Y_{\rho}}^2,
	\end{equation}
	with  $ \abs{u}_{Y_{\rho}}$ defined by \eqref{ynorm} and   $ \abs{u}_{X_{\rho}}$ in Definition \ref{defgev}.  
	
	For brevity of presentation, the proof of the statement \eqref{kees} is
 postponed to the end of this section.   We combine \eqref{kees} with the fact that
 \begin{eqnarray*}
 	\rho'^3\sum_{j=1}^2 \sum_{m=0}^{+\infty}L_{\rho,m}^2 \Big[   \frac{m+1}{\rho}\big( \norm{ \partial_t\partial_{x_j}^{m}u}_{L^2}^2+\norm{  \partial_y\partial_{x_j}^{m}u }_{L^2}^2 \big)+     \frac{(m+1)^3}{\rho^3}  \norm{\partial_{x_j}^{m}u}_{L^2}^2\Big] 
 	=\rho'^3\abs{u}_{Y_\rho}^2
 \end{eqnarray*}
 to conclude  for any $s\in[0,t]$ with the same time $t$ given in \eqref{has},  the following estimate holds 
 \begin{eqnarray*}
 	\begin{aligned}
 	&	 \rho'(s)^3\sum_{j=1}^2\sum_{m=0}^{+\infty}L_{\rho(s),m}^2 \Big[   \frac{m+1}{\rho(s)}\big( \norm{ \partial_t\partial_{x_j}^{m}u(s)}_{L^2}^2+\norm{  \partial_y\partial_{x_j}^{m}u(s) }_{L^2}^2 \big)+     \frac{(m+1)^3}{\rho(s)^3}  \norm{\partial_{x_j}^{m}u(s)}_{L^2}^2\Big]\\
 	&\quad +\sum_{j=1}^2\sum_{m=0}^{+\infty}L_{\rho(s),m}^2\big| \big( H_m(s), \    \partial_{x_j}^{m}u(s)\big)_{L^2} 	\big|+2\sum_{j=1}^2\sum_{m=0}^{+\infty}L_{\rho,m}^2\big|\big( H_m(s),\  \partial_t\partial_{x_j}^{m} u(s)\big)_{L^2}\big|\\
 	& \leq   \inner{  \rho'(s)^3+ C \rho(s)^{-2}   \abs{u(s)}_{X_{\rho(s)}}} \abs{u(s)}_{Y_{\rho(s)}}^2\leq - \Big( \frac{ \rho_0^3a^3}{8}  -32 \eps_0  C   \rho_0^{-2} \Big)e^{-s/32}  \abs{u(s)}_{Y_{\rho(s)}}^2\leq 0,
 	\end{aligned}
 \end{eqnarray*} 
 where the last inequality follows from the condition \eqref{has} and the property \eqref{prpof} by choosing  $\eps_0$  in \eqref{has} to be sufficiently small.   As a result, we take summation on both sides of \eqref{+pre}  for $j=1,2,$ 
   to obtain
 \begin{eqnarray*}
 	\begin{aligned}
		& \frac{1}{2} \frac{d}{dt} \sum_{j=1}^2\sum_{m=0}^{+\infty} L_{\rho,m}^2 \Big(2\norm{ \partial_t\partial_{x_j}^{m}u}_{L^2}^2+2\norm{  \partial_y\partial_{x_j}^{m}u }_{L^2}^2+\frac{d}{dt } \norm{ \partial_{x_j}^{m} u}_{L^2}^2+\norm{ \partial_{x_j}^{m} u}_{L^2}^2 \Big) \\
		&\quad+  \frac{d}{dt}   \sum_{j=1}^2\sum_{m=0}^{+\infty}  L_{\rho,m}^2\norm{   \partial_{x_j}^{m}u}_{L^2}^2 \Big(    \rho'^2 \frac{(m+1)^2}{\rho^2}-  \rho'\frac{m+1}{\rho} \Big)\\
		 &\quad+\frac18  \sum_{j=1}^2\sum_{m=0}^{+\infty}L_{\rho,m}^2\Big( \norm{  \partial_{x_j}^{m} u }_{L^2}^2+  \norm{ \partial_t\partial_{x_j}^{m} u}_{L^2}^2+ \norm{ \partial_y \partial_{x_j}^{m} u }_{L^2}^2 \Big)\\
		 &\quad  + 2 \sum_{j=1}^2\sum_{m=0}^{+\infty}  \rho'^2\frac{(m+1)^2}{\rho^2}  L_{\rho,m}^2 \norm{ \partial_{x_j}^{m} u}_{L^2}^2 \leq 0. 
	\end{aligned}
 \end{eqnarray*}
Thus, multiplying  the both sides by $e^{t/16}$ gives
 \begin{multline}\label{lete}
 	  \frac{1}{2} \frac{d}{dt}  e^{t/16} \sum_{j=1}^2\sum_{m=0}^{+\infty} L_{\rho,m}^2 \Big(2\norm{ \partial_t\partial_{x_j}^{m}u}_{L^2}^2+2\norm{  \partial_y\partial_{x_j}^{m}u }_{L^2}^2+\frac{d}{dt } \norm{ \partial_{x_j}^{m} u}_{L^2}^2+\norm{ \partial_{x_j}^{m} u}_{L^2}^2 \Big) \\
		 +  \frac{d}{dt} e^{t/16}  \sum_{j=1}^2 \sum_{m=0}^{+\infty} L_{\rho,m}^2\norm{   \partial_{x_j}^{m}u}_{L^2}^2  \Big(     \rho'^2 \frac{(m+1)^2}{\rho^2}-  \rho'\frac{m+1}{\rho} \Big) 
		 +\mathcal A\leq 0,
		 \end{multline}
		 with 
		 \begin{eqnarray*}
		 \begin{aligned}
		  \mathcal A&=-\frac{1}{32}e^{t/16} \sum_{j=1}^2\sum_{m=0}^{+\infty} L_{\rho,m}^2 \Big(2\norm{ \partial_t\partial_{x_j}^{m}u}_{L^2}^2+2\norm{  \partial_y\partial_{x_j}^{m}u }_{L^2}^2+\frac{d}{dt } \norm{ \partial_{x_j}^{m} u}_{L^2}^2+\norm{ \partial_{x_j}^{m} u}_{L^2}^2 \Big)\\
		 &\quad -\frac{1}{16}e^{t/16}   \sum_{j=1}^2\sum_{m=0}^{+\infty} L_{\rho,m}^2\norm{   \partial_{x_j}^{m}u}_{L^2}^2\Big(   \rho'^2 \frac{(m+1)^2}{\rho^2}-  \rho'\frac{m+1}{\rho} \Big)\\
		  &\quad +\frac18 e^{t/16} \sum_{j=1}^2\sum_{m=0}^{+\infty}L_{\rho,m}^2\Big( \norm{  \partial_{x_j}^{m} u }_{L^2}^2+  \norm{ \partial_t\partial_{x_j}^{m} u}_{L^2}^2+ \norm{ \partial_y \partial_{x_j}^{m} u }_{L^2}^2 \Big) \\
		  &\quad + 2 e^{t/16} \sum_{j=1}^2\sum_{m=0}^{+\infty}  \rho'^2\frac{(m+1)^2}{\rho^2}  L_{\rho,m}^2 \norm{ \partial_{x_j}^{m} u}_{L^2}^2. 
		  \end{aligned}
	\end{eqnarray*}
 By noting that $2|\rho'|(m+1)/\rho\leq 1+  \rho'^2 (m+1)^2/\rho^2$  and  
 \begin{multline}\label{equ}
 \frac14\inner{\norm{ \partial_t\partial_{x_j}^{m}u}_{L^2}^2+ \norm{  \partial_y\partial_{x_j}^{m}u }_{L^2}^2 +\norm{ \partial_{x_j}^{m} u}_{L^2}^2}\\ 
 \leq 	2\norm{ \partial_t\partial_{x_j}^{m}u}_{L^2}^2+2\norm{  \partial_y\partial_{x_j}^{m}u }_{L^2}^2+\frac{d}{dt } \norm{ \partial_{x_j}^{m} u}_{L^2}^2+\norm{ \partial_{x_j}^{m} u}_{L^2}^2\\
 \leq 3\inner{\norm{ \partial_t\partial_{x_j}^{m}u}_{L^2}^2+ \norm{  \partial_y\partial_{x_j}^{m}u }_{L^2}^2 +\norm{ \partial_{x_j}^{m} u}_{L^2}^2},
 \end{multline}
 we have 
 \begin{eqnarray*}
 	\mathcal A\geq 0.
 \end{eqnarray*}
 We now integrate  both sides of \eqref{lete} over $[0,s]$ for all $s\leq t$.   By the above inequality, \eqref{equ} and the fact that $\rho'\leq 0$,  we obtain
		\begin{equation}\label{ift}
		\begin{aligned}
		&  e^{ s/16}\sum_{j=1}^2\sum_{m=0}^{+\infty}L_{\rho(s) ,m}^2 \inner{ \norm{  \partial_{x_j}^{m} u(s) }_{L^2}^2+  \norm{ \partial_t\partial_{x_j}^{m} u(s) }_{L^2}^2 
		  + \norm{ \partial_y \partial_{x_j}^{m} u (s) }_{L^2}^2}   
		  \leq 12\big(|u_0|_{X_{\rho_0}}^2 +|u_1|_{X_{\rho_0}}^2\big)\\
		 &\quad +8\sum_{j=1}^2  \sum_{m=0}^{+\infty} L_{ \rho_0 ,m}^2\norm{   \partial_{x_j}^{m}u_0}_{L^2}^2 \Big(      \rho'(0)^2 \frac{ (m+1)^2}{ \rho_0^2}-\rho'(0)\frac{ m+1 }{ \rho_0} \Big)  \leq 16\big(|u_0|_{X_{2\rho_0}}^2 +|u_1|_{X_{2\rho_0}}^2\big),
		  \end{aligned}
		\end{equation}
  that is,  by Definition \ref{defgev} and \eqref{asonin},
	\begin{equation*}\label{inest}
	\forall\ s\in [0, t],\quad 	\abs{u(s)}_{X_{\rho(s)}} \leq  4e^{-s/32}\big(|u_0|_{X_{2\rho_0}}  +|u_1|_{X_{2\rho_0}} \big)\leq  4\eps_0e^{-s/32}.
	\end{equation*}
	Then
we have \eqref{last sta}.           
Hence,  $\mathbf C(t)$ in \eqref{assc} holds for all $t\geq 0$ by Proposition \ref{probootst} so that \eqref{arpe} holds.  We have proven the first assertion \eqref{arpe} in Theorem \ref{thmofa}.

 Now it remains to prove \eqref{kees}           as follows.
 
 \begin{proof}[\bf Proof of   \eqref{kees}]
    Recall $H_m$  given by \eqref{eqmu}. 
 	We first estimate the terms  involving the pressure function. Firstly, note that  the following compatibility condition for the solution $u$ to \eqref{hyper} holds
 \begin{equation}
 	\label{cocon}
 	\forall\ x\in\mathbb R^2,\quad \int_0^1\partial_x\cdot u (x,y)dy=0.
 \end{equation}
  Then integrating by parts gives
	\begin{equation}\label{dfc}
		\big|\big( \partial_x\partial_{x_j}^mp,\ \partial_{x_j}^{m} u \big)_{L^2}\big|=\big|\big( \partial_{x_j}^mp,\ \partial_{x_j}^{m}\partial_x \cdot u \big)_{L^2}\big|=\big|\big( \partial_{x_j}^mp,\ \partial_y \partial_{x_j}^{m} v\big)_{L^2}\big|=0, 
	\end{equation}
	where the last equality follows from integration by parts and the  fact that $\partial_y p=0$.  
 The above equality also holds  when $\partial_{x_j}^m u$ is replaced by $\partial_t\partial_{x_j}^m u$. 
 That is, 
 \begin{equation}\label{pres}
 	 \big|\big( \partial_x\partial_{x_j}^mp,\ \partial_{x_j}^{m} u \big)_{L^2}\big|+ \big|\big( \partial_x\partial_{x_j}^mp,\  \partial_t  \partial_{x_j}^{m} u \big)_{L^2}\big|=0.
 \end{equation}
 It remains to estimate the other terms in $ H_m$ and show that 
 \begin{equation}\label{fines}
 \begin{aligned}
		&\sum_{m=0}^{+\infty}L_{\rho,m}^2\Big|\Big( \sum_{k=0}^{m}\frac{m!}{k!(m-k)!}  \big[(\partial_{x_j}^k u  \cdot\partial_{x})\partial_{x_j}^{m-k} u+(\partial_{x_j}^k v)\partial_{x_j}^{m-k}\partial_y u \big],\ \partial_t \partial_{x_j}^{m} u \Big)_{L^2}\Big|\\
		&+\sum_{m=0}^{+\infty}L_{\rho,m}^2\Big|\Big( \sum_{k=0}^{m}\frac{m!}{k!(m-k)!}  \big[(\partial_{x_j}^k u  \cdot\partial_{x})\partial_{x_j}^{m-k} u+(\partial_{x_j}^k v)\partial_{x_j}^{m-k}\partial_y u \big],\   \partial_{x_j}^{m} u \Big)_{L^2}\Big| \leq  C     \abs{u}_{X_{\rho}}\abs{u}_{Y_{\rho}}^2.
		\end{aligned}
	\end{equation}
 To prove \eqref{fines}, we  write
 \begin{equation}\label{s12}
	\sum_{m=0}^{+\infty}L_{\rho,m}^2\big|\big( \sum_{k=0}^{m}\frac{m!}{k!(m-k)!}  \big[(\partial_{x_j}^k u  \cdot\partial_{x})\partial_{x_j}^{m-k} u+(\partial_{x_j}^k v)\partial_{x_j}^{m-k}\partial_y u \big],\ \partial_t \partial_{x_j}^{m} u \big)_{L^2}\big| \leq S_1+S_2,
\end{equation} 
where
\begin{equation}\label{s1s2}
\left\{
	\begin{aligned}
S_1&=\sum_{m=0}^{+\infty} \sum_{k=0}^{m}  L_{\rho,m}^2 \frac{m!}{k!(m-k)!}      \norm{   (\partial_{x_j}^k u\cdot\partial_x) \partial_{x_j}^{m-k}u}_{L^2}\norm{ \partial_t \partial_{x_j}^{m} u}_{L^2},\\  
	S_2&= \sum_{m=0}^{+\infty}   \sum_{k=0}^{m}   L_{\rho,m}^2\frac{m!}{k!(m-k)!}  \norm{   (\partial_{x_j}^k v) \partial_y\partial_{x_j}^{m-k}u}_{L^2}\norm{ \partial_t \partial_{x_j}^{m} u}_{L^2}. 
	\end{aligned}
\right.
\end{equation} 
We first estimate  $S_1$ as follows.
\begin{equation}\label{s1}
	\begin{aligned}
		S_1 &\leq \sum_{m=0}^{+\infty} \sum_{k=0}^{[m/2]}  L_{\rho,m}^2 \frac{m!}{k!(m-k)!}      \norm{  \partial_{x_j}^k u}_{L^\infty}\norm{   \partial_x \partial_{x_j}^{m-k}u}_{L^2}\norm{ \partial_t \partial_{x_j}^{m} u}_{L^2}\\
		&\quad+\sum_{m=0}^{+\infty} \sum_{k=[m/2]+1}^{m}  L_{\rho,m}^2 \frac{m!}{k!(m-k)!}      \norm{    \partial_{x_j}^k u  }_{L^2}\norm{ \partial_x \partial_{x_j}^{m-k}u}_{L^\infty}\norm{ \partial_t \partial_{x_j}^{m} u}_{L^2}\\
		&:=S_{1,1}+S_{1,2},
	\end{aligned}
\end{equation}
   where   $[m/2] $ stands for  the largest integer   $\leq m/2.$  
 To estimate $S_{1,1}$	 and  $S_{1,2}$,  we will use the following inequalities that follow from 
 straightforward calculation.
 If  $0\leq k\leq [m/2]$, then  
 \begin{equation}\label{fac1}
 		\begin{aligned}
 			   & \frac{m! }{ k!(m-k)!} \frac{ L_{\rho,m}  }{ L_{\rho,k+2} L_{\rho,m-k+1}}   = \frac{m!} {k!(m-k)!} \frac{[(k+2)!]^2}{\rho^{k+3}(k+3)^7} \frac{[(m-k+1)!]^2}{  \rho^{m-k+2} (m-k+2)^7} \frac{ \rho^{m+1}(m+1)^7}{(m!)^2 }  \\
 			     &\leq C \frac{k! (m-k+1)! (m-k+1)}{m!(k+3)\rho^{4}}     \leq C  \frac{ (m-k+1)m } { (k+1)\rho^{4} }    \leq \frac{C}{ \rho^2} \frac{1}{k+1 }\frac{(m+1)^{1/2}}{\rho^{1/2}}\frac{(m-k+2)^{3/2}}{\rho^{3/2}}.
 			      		\end{aligned}
 	\end{equation}
 On the other hand,  if $  [m/2]+1\leq k\leq m$, then
\begin{equation}\label{fac2}
 			    \frac{m! }{ k!(m-k)!} \frac{ L_{\rho,m}  }{  L_{\rho,k} L_{\rho,m-k+3}}          \leq   \frac{C} { \rho^4}\frac{1} {  m-k+1}.   
 \end{equation}
Recalling $S_{1,1}$  given in \eqref{s1},
 by \eqref{fac1} and  the definition of $\abs{\cdot}_{Y_\rho}$  in \eqref{ynorm}, we have
 \begin{multline}\label{fre}
 \begin{aligned}
 	&S_{1,1} =  \sum_{m=0}^{+\infty} \sum_{k=0}^{[m/2]}   \frac{m!}{k!(m-k)!}\frac{L_{\rho,m}}{ L_{\rho,k+2}L_{\rho,m-k+1}}      \Big(L_{\rho,k+2}\norm{  \partial_{x_j}^k u}_{L^\infty}\Big)   \Big(L_{\rho,m-k+1}\norm{    \partial_x\partial_{x_j}^{m-k}u}_{L^2}\Big)\\
 	&\qquad\qquad \qquad\quad\times L_{\rho,m} \norm{ \partial_t \partial_{x_j}^{m} u}_{L^2}\\
 	&\leq   \frac{C}{\rho^2}\sum_{m=0}^{+\infty} \sum_{k=0}^{[m/2]}         \frac{L_{\rho, k+2}\norm{  \partial_{x_j}^k u}_{L^\infty}}{k+1} \\
 	&\qquad\quad \times\Big( \frac{(m-k+2)^{3/2}}{\rho^{3/2}}L_{\rho,m-k+1}\norm{    \partial_x\partial_{x_j}^{m-k}u}_{L^2} \Big) \Big(\frac{(m+1)^{1/2}}{\rho^{1/2}}L_{\rho,m} \norm{ \partial_t \partial_{x_j}^{m} u}_{L^2}\Big)\\
 	&\leq  \frac{ C}{\rho^2} \bigg[\sum_{m=0}^{+\infty} \Big(\sum_{k=0}^{[m/2]}         \frac{L_{\rho,k+2}\norm{  \partial_{x_j}^k u}_{L^\infty}}{k+1} \times   \frac{(m-k+2)^{3/2}}{\rho^{3/2}} L_{\rho,m-k+1}\norm{\partial_x  \partial_{x_j}^{m-k}u}_{L^2}  \Big)^2\bigg]^{1\over 2}\abs{u}_{Y_{\rho}}\\
 	&\leq \frac{ C}{\rho^2}   \sum_{k=0}^{\infty}         \frac{L_{\rho,k+2}\norm{  \partial_{x_j}^k u}_{L^\infty}}{k+1}  \bigg[\sum_{m=0}^{+\infty}     \frac{ (m+2)^3}{\rho^3}  L_{\rho, m+1}^2\norm{ \partial_x  \partial_{x_j}^{ m}u}_{L^2}^2  \bigg]^{1/2}  \abs{u}_{Y_{\rho}},
\end{aligned}
 \end{multline}
 where we have used   Young's inequality for   discrete convolution. Moreover, it follows from \eqref{equofn} that  
  \begin{eqnarray*}
 	\bigg[\sum_{m=0}^{+\infty}     \frac{(m+2)^3}{\rho^3}  L_{\rho, m+1}^2\norm{ \partial_x  \partial_{x_j}^{m}u}_{L^2}^2  \bigg]^{1\over 2} 
 	  \leq  \bigg[\sum_{m=0}^{+\infty}    \frac{ (m+1)^3}{\rho^3}  L_{ \rho, m }^2\inner{\norm{  \partial_{x_1}^{m}u}_{L^2}^2+\norm{  \partial_{x_2}^{m}u}_{L^2}^2}  \bigg]^{1\over 2} \leq  \abs{u}_{Y_\rho}. \end{eqnarray*}
  By the Sobolev embedding inequality
 \begin{equation*}
	 \norm{F}_{L^\infty}\leq  C \big(\norm{F}_{H_x^2(L_y^2)}+\norm{\partial_y F}_{H_x^2(L_y^2)} \big) 
	\end{equation*}
  and \eqref{equofn},  it holds that
 	\begin{multline*}
 		 \sum_{k=0}^{\infty}         \frac{L_{\rho,k+2}\norm{  \partial_{x_j}^k u}_{L^\infty}}{k+1}\leq C \sum_{k=0}^{\infty}         \frac{1}{k+1}L_{\rho,k+2}\inner{\norm{  \partial_{x_j}^k u}_{H_x^2(L_y^2)}+\norm{  \partial_y\partial_{x_j}^{k} u}_{H_x^2(L_y^2)} }\\
 		 \leq C \bigg[\sum_{k=0}^{\infty}          L_{\rho,k+2}^2\big(\norm{  \partial_{x_j}^k u}_{H_x^2(L_y^2)}^2+\norm{  \partial_y\partial_{x_j}^{k} u}_{H_x^2(L_y^2)}^2 \big) \bigg]^{1/2}\leq C \abs{u}_{X_{\rho}}.
 	\end{multline*} 
 Hence, combining the above estimates with \eqref{fre} yields
 \begin{equation*}\label{stcom}
 	S_{1,1} \leq    C  \rho^{-2}   \abs{u}_{X_{\rho}}\abs{u}_{Y_{\rho}}^2.
 \end{equation*}
 	Similarly,   
 	by using \eqref{fac2}, we have
 	\begin{eqnarray*}
 		S_{1,2} \leq C \rho^{-2}  \abs{u}_{X_{\rho}}\abs{u}_{Y_{\rho}}^2.
 	\end{eqnarray*}
 Substituting the  estimates on $S_{1,1}$ and $S_{1,2}$ into \eqref{s1} yields
 \begin{equation*}
 	S_1\leq C \rho^{-2} \abs{u}_{X_{\rho}}\abs{u}_{Y_{\rho}}^2.
 \end{equation*}
The term $S_2$  in \eqref{s1s2} can be estimated similarly and we
omit the details for brevity, that is,  
  \begin{equation*}
 	S_2\leq C\rho^{-2} \abs{u}_{X_{\rho}}\abs{u}_{Y_{\rho}}^2. 
 \end{equation*}
 In summary, we have, in view of  \eqref{s12}, 
 \begin{eqnarray*}
 	\sum_{m=0}^{+\infty}L_{\rho,m}^2\big|\big( \sum_{k=0}^{m}\frac{m!}{k!(m-k)!}  \big[(\partial_{x_j}^k u  \cdot\partial_{x})\partial_{x_j}^{m-k} u+(\partial_{x_j}^k v)\partial_{x_j}^{m-k}\partial_y u \big], \partial_t \partial_{x_j}^{m} u \big)_{L^2}\big| \leq C\rho^{-2} \abs{u}_{X_{\rho}}\abs{u}_{Y_{\rho}}^2.
 \end{eqnarray*}
 Similarly, 
 \begin{eqnarray*}
 	\sum_{m=0}^{+\infty}L_{\rho,m}^2\big|\big( \sum_{k=0}^{m}\frac{m!}{k!(m-k)!}  \big[(\partial_{x_j}^k u  \cdot\partial_{x})\partial_{x_j}^{m-k} u+(\partial_{x_j}^k v)\partial_{x_j}^{m-k}\partial_y u \big],\  \partial_{x_j}^{m} u \big)_{L^2}\big| \leq C\rho^{-2} \abs{u}_{X_{\rho}}\abs{u}_{Y_{\rho}}^2.
 \end{eqnarray*}
Therefore,  the statement \eqref{kees} holds.   
\end{proof}

\subsection{Proof of Theorem \ref{thmofa}: The second assertion}  
Since the argument is similar to the one used
  in the previous section,  we now sketch the proof of 
  \eqref{secass} in Theorem \ref{thmofa}  for brevity.    In fact, using the notation 
\begin{eqnarray*}
	\mathcal U=(u,\partial_tu),  \quad \mathcal H_m=(H_m, \partial_tH_m)
\end{eqnarray*}
with $H_m$  defined in \eqref{eqmu},  we have 
by applying $\partial_t\partial_{x_j}^m $ to \eqref{hyper}, 
\begin{eqnarray}\label{eqmac}
	 \begin{aligned} 
	  \big(\partial_t^2+\partial_t  -\partial_y^2\big)   \partial_{x_j}^{m} \mathcal U&=\mathcal H_m.
	 \end{aligned}
\end{eqnarray}
Similar to \eqref{kees}, we conclude 
\begin{equation}\label{laet}
		\sum_{j=1}^2\sum_{m=0}^{+\infty}L_{\rho,m}^2\big| \big( \mathcal H_m, \     \partial_{x_j}^{m}\mathcal U\big)_{L^2} 	\big|+2\sum_{j=1}^2\sum_{m=0}^{+\infty}L_{\rho,m}^2\big|\big(  \mathcal H_m,\  \partial_t\partial_{x_j}^{m} \mathcal U\big)_{L^2}\big| \leq  C   \rho^{-2}   \abs{\mathcal U}_{X_{\rho}}\abs{\mathcal U}_{Y_{\rho}}^2.
	\end{equation}
 In fact, as  \eqref{dfc},  the  equations 
\begin{eqnarray*}
\big( \partial_t\partial_x\partial_{x_j}^mp,\ \partial_{x_j}^{m} u \big)_{L^2} =\big( \partial_t\partial_x\partial_{x_j}^mp,\ \partial_t\partial_{x_j}^{m} u \big)_{L^2} = \big( \partial_t\partial_x\partial_{x_j}^mp,\ \partial_t^2   \partial_{x_j}^{m} u\big)_{L^2}=0
\end{eqnarray*}
also hold. Hence, 
\begin{eqnarray*}
	\big( \partial_t\partial_x\partial_{x_j}^mp,\ \partial_{x_j}^{m}\mathcal U +2\partial_t\partial_{x_j}^{m}\mathcal U \big)_{L^2}=0.
\end{eqnarray*}
As a result, \eqref{laet} follows by repeating the  argument for obtaining \eqref{kees}.  

 We now
take the inner product with $ \partial_{x_j}^m\mathcal U+2\partial_t\partial_{x_j}^m\mathcal U$  on both sides of the   equation \eqref{eqmac}, and then repeat the argument for proving  \eqref{arpe} with $\partial_{x_j}^mu$   replaced by $ \partial_{x_j}^m\mathcal U$. By   \eqref{laet},
we have   the following  estimate that is  similar to \eqref{ift}:
	\begin{equation}\label{lut}
		\begin{aligned}
	e^{ t/16}\abs{\mathcal U(t)}_{X_\rho}^2	&=  e^{ t/16}\sum_{j=1}^2\sum_{m=0}^{+\infty}L_{\rho(t) ,m}^2 \inner{ \norm{  \partial_{x_j}^{m} \mathcal U(t) }_{L^2}^2+  \norm{ \partial_t\partial_{x_j}^{m}\mathcal U(t) }_{L^2}^2 
		  + \norm{ \partial_y \partial_{x_j}^{m}\mathcal U (t) }_{L^2}^2}    \\
		  & \leq 12 \sum_{j=1}^2\sum_{m=0}^{+\infty}L_{\rho_0 ,m}^2 \inner{ \norm{  \partial_{x_j}^{m} \mathcal U|_{t=0} }_{L^2}^2+  \norm{ \partial_t\partial_{x_j}^{m}\mathcal U|_{t=0} }_{L^2}^2 
		  + \norm{ \partial_y \partial_{x_j}^{m}\mathcal U|_{t=0} }_{L^2}^2} \Big]\\
		 &\quad +8\sum_{j=1}^2  \sum_{m=0}^{+\infty} L_{ \rho_0 ,m}^2\norm{   \partial_{x_j}^{m}\mathcal U|_{t=0}}_{L^2}^2 \Big(      \rho'(0)^2 \frac{ (m+1)^2}{ \rho_0^2}-\rho'(0)\frac{ m+1 }{ \rho_0} \Big) \\
		 & \leq C \big(|u_0|_{X_{2\rho_0}}^2 +|u_1|_{X_{2\rho_0}}^2\big)+12\sum_{j=1}^2 \sum_{m=0}^{+\infty}L_{\rho_0 ,m}^2     \norm{ \partial_t^2\partial_{x_j}^{m}u|_{t=0} }_{L^2}^2,
		  \end{aligned}
		\end{equation}
	where we have used the fact in	the last inequality 
		\begin{eqnarray*}
	\mathcal U|_{t=0}=( u_0, u_1),\quad		\partial_y \mathcal U|_{t=0}=(\partial_yu_0,\partial_yu_1),\quad \partial_t  \mathcal U|_{t=0}=( u_1, \partial_t^2u|_{t=0}).
		\end{eqnarray*}
It remains to estimate the last term on the right hand side of \eqref{lut}. For this,
 when $t=0$,  the first equation of \eqref{hyper} gives
\begin{equation}\label{trac}
	\partial_t^2u|_{t=0}=-u_1-(u_0\cdot\partial_x)u_0 -v_0\partial_yu_0+\partial_y^2u_0-\partial_xp|_{t=0}, \quad v_0=-\int_0^y \partial_x\cdot u_0\,d\tilde y.
\end{equation}
By applying the divergence operator  to the first equation in \eqref{hyper} and then using the compatibility condition \eqref{cocon},   we have the  elliptic equation for the pressure function
\begin{equation} \label{ellpre}
	\Delta_xp=-  \partial_x\cdot   \int_0^1  \big[(u\cdot \partial_x) u +  ( \partial_x\cdot u) u \big] dy   +\partial_x\cdot \int_0^1 \partial_y^2 u dy, \quad x\in\mathbb R^2.
\end{equation} 
Thus,  the   standard elliptic theory implies that 
\begin{equation}\label{eforpressure}
	\sum_{m=0}^\infty L_{\rho_0, m}^2 \norm{\partial_{x_j}^m\partial_x p|_{t=0}}_{L_x^2}^2\leq C\big(\abs{u_0}_{X_{2\rho_0}}^4+\abs{\partial_yu_0}_{X_{2\rho_0}}^2\big),
\end{equation}
This  together with \eqref{trac} yields
\begin{equation}\label{trau}
	\sum_{m=0}^\infty L_{\rho_0, m}^2 \norm{\partial_t^2\partial_{x_j}^m u|_{t=0}}_{L^2}^2\leq C\big(\abs{u_0}_{X_{2\rho_0}}^4+\abs{\partial_yu_0}_{X_{2\rho_0}}^2+\abs{ u_1}_{X_{2\rho_0}}^2\big).
\end{equation}
For completeness, we give  the proof of \eqref{eforpressure} and \eqref{trau} in Appendix \ref{secap}. 
Combining the above estimate with \eqref{lut} gives the second assertion in Theorem  \ref{thmofa}.

\subsection{Proof of Theorem \ref{thmofa}: The third assertion} It remains to prove \eqref{yyu}. 
  We take the $L^2$ inner product with 
 $\partial_{x_j}^m \partial_y^2u$   on both sides of the equation  
 \begin{equation*}	
 \begin{aligned} 
	   \big(\partial_t^2+\partial_t  -\partial_y^2\big)  \partial_{x_j}^{m} u&=-\partial_x\partial_{x_j}^mp -\underbrace{ \sum_{k=0}^{m}{m\choose k} \Big[(\partial_{x_j}^k u  \cdot\partial_{x})\partial_{x_j}^{m-k} u+(\partial_{x_j}^k v)\partial_{x_j}^{m-k}\partial_y u \Big]}_{\stackrel{def}{ =}Q_m},  
	 \end{aligned}
\end{equation*}
and then multiply  by $L_{\rho/2,m}^2$ before taking  summation for $m$.  This gives
\begin{multline*}
   \sum_{m=0}^{+\infty} L_{\rho/2,m}^2\norm{ \partial_y^2 \partial_{x_j}^{m} u }_{L^2}^2  \leq  \sum_{m=0}^{+\infty} L_{\rho/2,m}^2\big(\norm{ \partial_t^2 \partial_{x_j}^{m} u }_{L^2}+\norm{ \partial_t \partial_{x_j}^{m} u }_{L^2}\big)\norm{ \partial_y^2 \partial_{x_j}^{m} u }_{L^2}   \\
+  \sum_{m=0}^{+\infty} L_{\rho/2,m}^2 \norm{Q_m}_{L^2} \norm{\partial_y^2\partial_{x_j}^{m}u}_{L^2}   
	 + \sum_{m=0}^{+\infty} L_{\rho/2,m}^2  \norm{\partial_{x_j}^{m}\partial_xp}_{L_x^2} \norm{ \partial_y^2 \partial_{x_j}^{m} u }_{L^2}.
\end{multline*}
Thus
 \begin{multline}\label{eqyu}
	  \sum_{m=0}^{+\infty} L_{\rho/2,m}^2\norm{ \partial_y^2 \partial_{x_j}^{m} u }_{L^2}^2   \leq C  \big( \abs{\partial_tu}_{X_{\rho/2}}^2+\abs{u}_{X_{\rho/2}}^2\big)   \\
+4  \sum_{m=0}^{+\infty} L_{\rho/2,m}^2 \norm{Q_m}_{L^2} \norm{\partial_y^2\partial_{x_j}^{m}u}_{L^2}  
	 +2 \sum_{m=0}^{+\infty} L_{\rho/2,m}^2  \norm{\partial_{x_j}^{m}\partial_xp}_{L_x^2}^2.
\end{multline} 
 Using a similar argument as in  \eqref{appco} of Appendix \ref{secap}  gives
\begin{multline*}
 \sum_{m=0}^{+\infty} L_{\rho/2,m}^2 \norm{Q_m}_{L^2} \norm{\partial_y^2\partial_{x_j}^{m}u}_{L^2}  \leq C    \abs{u}_{X_\rho}^2 \Big(\sum_{m=0}^{+\infty}   L_{\rho/2,m}^2\norm{ \partial_y^2 \partial_{x_j}^{m} u }_{L^2}^2\Big)^{1\over 2} \\
  \leq \frac{1}{8}      \sum_{m=0}^{+\infty}   L_{\rho/2,m}^2\norm{ \partial_y^2 \partial_{x_j}^{m} u }_{L^2}^2  +  C    \abs{u}_{X_\rho}^4.
  \end{multline*}
We use the elliptic theory for \eqref{ellpre} and a similar computation as in \eqref{eforpressure}  to conclude
\begin{multline*}
	  \sum_{m=0}^{+\infty} L_{\rho/2,m}^2   \norm{\partial_{x_j}^{m}\partial_xp}_{L_x^2}^2 
	 \leq   C  \abs{u}_{X_\rho}^4 +C   \sum_{m=0}^{+\infty} L_{\rho/2,m}^2 \norm{\partial_y\partial_{x_j}^{m} u}_{L_x^2L_y^\infty}^2 \\
	 \leq  \frac18  \sum_{m=0}^{+\infty} L_{\rho/2,m}^2  \norm{ \partial_y^2 \partial_{x_j}^{m} u }_{L^2}^2+C   \big(\abs{u}_{X_\rho}^4+  \abs{u}_{X_\rho}^2\big), 
\end{multline*}
Here, in the last inequality, we have used  the fact that 
\begin{eqnarray*}
	\forall\ r\in[0,1],\quad   \big(\partial_yu(r))^2=\int_\xi^r \partial_y  \big(\partial_yu( y))^2d  y=2\int_\xi^r  \big(\partial_y^2u(  y) \big)\partial_y   u( y) d  y,
\end{eqnarray*}
with $\xi\in[0,1]$ satisfying 
\begin{eqnarray*}
	\partial_yu(\xi)=\int_0^1\partial_y u dy=0
\end{eqnarray*}
because of the boundary condition $u|_{y=0,1}=0$.   Substituting the above estimates into \eqref{eqyu} yields
 \begin{eqnarray*}
   \sum_{m=0}^{+\infty} L_{\rho/2,m}^2\norm{ \partial_y^2 \partial_{x_j}^{m} u }_{L^2}^2\leq C   \big(\abs{u}_{X_\rho}^2+  \abs{u}_{X_\rho}^4+ \abs{\partial_tu}_{X_\rho}^2\big).
 \end{eqnarray*}
 Observe
 \begin{eqnarray*}
 	 \abs{\partial_yu}_{X_{\rho/2}}^2\leq   \sum_{m=0}^{+\infty} L_{\rho/2,m}^2\norm{ \partial_y^2 \partial_{x_j}^{m} u }_{L^2}^2+ \abs{u}_{X_{\rho/2}}^2+ \abs{\partial_tu}_{X_{\rho/2}}^2,
 \end{eqnarray*}
so that   \eqref{yyu} follows.   The proof of  Theorem \ref{thmofa}  
is completed.

   \section{ Global well-posedness of   original system}
   
   We will prove Theorem  \ref{thmhas} about the global well-posedness of the anisotropic hyperbolic Navier-Stokes system \eqref{hans} in this section.  The argument is  similar to the one used in  the previous section so that we only give the sketch of the proof.
  
\begin{proof}
[{\bf Sketch of the proof of Theorem \ref{thmhas}}]   	 We apply $\partial_{x_j}$ to the evolution equations of $u^\eps$ and $v^\eps$  in \eqref{hans} to obtain
    \begin{multline} \label{tme}
	  \big(\partial_t^2+\partial_t -\eps^2 \Delta_x -\partial_y^2\big)  \partial_{x_j}^{m} u^\eps \\
	  =-\partial_x\partial_{x_j}^mp^\eps -\sum_{k=0}^{m}{m\choose k} \Big[(\partial_{x_j}^k u^\eps  \cdot\partial_{x})\partial_{x_j}^{m-k} u^\eps+(\partial_{x_j}^k v^\eps)\partial_{x_j}^{m-k}\partial_y u^\eps \Big]  
	  \stackrel{def}{=} T_m^\eps, 
	 \end{multline}
    and
    \begin{multline} \label{nme}
	 \eps^2 \big(\partial_t^2+\partial_t -\eps^2 \Delta_x -\partial_y^2\big)  \partial_{x_j}^{m}  v^\eps\\
	 =-\partial_y\partial_{x_j}^mp^\eps -\eps^2\sum_{k=0}^{m}{m\choose k} \Big[(\partial_{x_j}^k u^\eps  \cdot\partial_{x})\partial_{x_j}^{m-k} v^\eps+(\partial_{x_j}^k v^\eps)\partial_{x_j}^{m-k}\partial_y v^\eps \Big]  
	  \stackrel{def}{=} N_m^\eps. 
	 \end{multline}
   It follows from the divergence-free condition that
     \begin{eqnarray*}
  	  \big( \partial_x\partial_{x_j}^mp^\eps,\ \partial_{x_j}^{m} u^\eps \big)_{L^2}+\big( \partial_y\partial_{x_j}^mp^\eps,\ \partial_{x_j}^{m}v^\eps \big)_{L^2}=\big( \partial_x\partial_{x_j}^mp^\eps, \partial_t\partial_{x_j}^{m} u^\eps \big)_{L^2}+\big( \partial_y\partial_{x_j}^mp^\eps, \partial_t\partial_{x_j}^{m}v^\eps \big)_{L^2}=0.
  \end{eqnarray*} 
   Let $T^\eps_m, N^\eps_m$ be  in \eqref{tme} and \eqref{nme}. 
  Then by using the above equality  instead of \eqref{pres} and  following the argument in the proof of  \eqref{kees}, we conclude
    \begin{multline*}
    		\sum_{j=1}^2\sum_{m=0}^{+\infty}L_{\rho,m}^2  \big( T_m^\eps, \    \partial_{x_j}^{m}u^\eps+2\partial_t\partial_{x_j}^{m} u^\eps\big)_{L^2} 	   
    		+\sum_{j=1}^2\sum_{m=0}^{+\infty}L_{\rho,m}^2  \big( N_m^\eps,   \partial_{x_j}^{m}v^\eps+2 \partial_t\partial_{x_j}^{m} v^\eps\big)_{L^2} 	 \\
    		 \leq  C   \rho^{-2}   \abs{(u^\eps, \eps v^\eps)}_{X_{\rho}}\abs{(u^\eps, \eps v^\eps)}_{Y_{\rho}}^2,
    \end{multline*}
    where $\abs{(u^\eps, \eps v^\eps)}_{X_{\rho}}$  and $\abs{(u^\eps, \eps v^\eps)}_{Y_{\rho}}$ are given by Definition \ref{defgev} and \eqref{ynorm}.
     Then we 
     take the $L^2$ inner product with $\partial_{x_j}^m u^\eps+2 \partial_t\partial_{x_j}^m u^\eps$ in  \eqref{tme}  and with $\partial_{x_j}^m v^\eps+2 \partial_t\partial_{x_j}^m v^\eps$ in \eqref{nme}.
  The following  
a priori estimate for the system \eqref{hans} can be obtained by using the argument in the previous
section:
\begin{equation*}
	\forall\ t\geq 0,\quad 	\abs{\big(u^\eps(t), \eps v^\eps(t)\big) }_{X_{\rho(t)}} \leq  4\delta_0e^{-t/32},
	\end{equation*}
provided that the initial data satisfy 
  $$
 	|(u_0^\eps, \eps v_0^\eps)|_{X_{2\rho_0}}  +|(u_1^\eps, \eps v_1^\eps)|_{X_{2\rho_0}} \leq \delta_0
$$ with $\delta_0$  being sufficiently small. Hence,  the proof of  Theorem  \ref{thmhas} is completed.

\end{proof}

\section{Hydrostatic limit }

In the final  section,  we will prove Theorem \ref{thmlim} about  the hydrostatic limit from   \eqref{hans} to \eqref{hyper} as $\eps\rightarrow 0$.    

Suppose the assumptions in  Theorem  \ref{thmlim}  hold. Let    $(u^\eps, v^\eps, p^\eps) $ and $(u,p)$  solve the  anisotropic hyperbolic  Navier-Stokes system \eqref{hans} and the hyperbolic hydrostatic Navier-Stokes system \eqref{hyper}, respectively. In the following discussion,     
we use the notation 
 \begin{equation*}
	U^\eps=u^\eps-u,\quad V^\eps=  v^\eps-v,\quad P^\eps=p^\eps-p,
\end{equation*}
where $v(t,x,y)=-\int_0^y\partial_x\cdot u (t,x,\tilde y)d\tilde y$.  Then it follows from \eqref{hans} and \eqref{hyper} that  \begin{equation*}
	\left\{
\begin{aligned}
&  \big(\partial_t^2+\partial_t -\eps^2 \Delta_x -\partial_y^2\big) U^\eps + \partial_x P^\eps
=\eps^2\Delta_x u+R_{\eps},    \\
& \eps^2 \big(\partial_t^2+\partial_t     -\eps^2\Delta_x-\partial_y^2\big) V ^\eps + \partial_yP^\eps=-\eps^2 \big(\partial_t^2+\partial_t +u\cdot\partial_x+v\partial_y    -\eps^2\Delta_x-\partial_y^2\big) v+S_\eps, \\
&\partial_x\cdot U^\eps +\partial_yV^\eps =0, \\
&U^\eps|_{y=0,1}=0,\quad V^\eps|_{y=0,1}=0,\\
&(U^\eps, V^\eps)|_{t=0}=(u^\eps_0-u_0, v^\eps_0-v_0),\quad (\partial_tU^\eps, \partial_tV^\eps)|_{t=0}=(u^\eps_1-u_1, v^\eps_1-v_1),
\end{aligned}\right.
\end{equation*}
where $v_j=-\int_0^y\partial_x\cdot u_j (t,x,\tilde y)d\tilde y,  j=0,1,$ and 
\begin{eqnarray*}
\begin{aligned}
	R_{\eps}&= -(U^\eps\cdot\partial_x)u^\eps- (u\cdot\partial_x)U^\eps  -V^\eps\partial_yu^\eps-v\partial_yU^\eps,\\
	S_\eps &=-\eps^2(U^\eps\cdot\partial_x)v^\eps-\eps^2(u\cdot\partial_x)V^\eps-\eps^2V^\eps\partial_yv^\eps-\eps^2v\partial_yV^\eps.  
		\end{aligned}
\end{eqnarray*}
In the following, we will use $C$ to denote a
generic constant depending only on $\rho_0$ and the Sobolev embedding constant but independent of $\eps$.   
Similar to the argument used in Section \ref{secderivation}, we have
 \begin{equation}\label{st1}
\sum_{m=0}^{+\infty}L_{\rho/2,m}^2	\big(\partial_{x_j}^m R_\eps,\   \partial_{x_j}^m U^\eps+2\partial_t\partial_{x_j}^m U^\eps \big)_{L^2}\leq C\inner{\abs{u^\eps}_{X_\rho}+\abs{u}_{X_\rho}}\abs{U^\eps}_{Y_{\rho/2}}^2,
\end{equation}
and
\begin{multline}\label{st2}
\sum_{m=0}^{+\infty}L_{\rho/2,m}^2	\big(\partial_{x_j}^m S_\eps,\   \partial_{x_j}^m V^\eps+2\partial_t\partial_{x_j}^m V^\eps \big)_{L^2}\\ 
\leq C \abs{\eps v^\eps}_{X_\rho} \abs{U^\eps}_{Y_{\rho/2}}^2 +C\inner{\abs{\eps v^\eps}_{X_\rho}+\abs{u^\eps}_{X_\rho}+\abs{u}_{X_\rho}}\abs{\eps V^\eps}_{Y_{\rho/2}}^2,
\end{multline}
where $\abs{\cdot}_{Y_{\rho/2}}$ is defined in \eqref{ynorm}. Please refer to   Appendix \ref{secap} for the detailed proof.     
 Moreover, direct calculation gives
 \begin{eqnarray*}
 	\sum_{m=0}^{+\infty}L_{\rho/2,m}^2	\big(\eps^2\partial_{x_j}^m\Delta_x u,\   \partial_{x_j}^m U^\eps+2\partial_t\partial_{x_j}^m U^\eps \big)_{L^2}\leq C\eps^2 \abs{u}_{X_\rho} \abs{U^\eps}_{X_{\rho/2}},
 \end{eqnarray*}
 and
 \begin{multline}\label{vsot}
 	\sum_{m=0}^{+\infty}L_{\rho/2,m}^2	\big(-\eps^2\partial_{x_j}^m \big(\partial_t^2+\partial_t +u\cdot\partial_x+v\partial_y    -\eps^2\Delta_x-\partial_y^2\big) v,\   \partial_{x_j}^m V^\eps+2\partial_t\partial_{x_j}^m V^\eps \big)_{L^2}\\
 	\leq C\eps \inner{ \abs{u}_{X_\rho} + \abs{u}_{X_\rho}^2+ \abs{\partial_tu}_{X_{\rho}}} \abs{\eps V^\eps}_{X_{\rho/2}},
 \end{multline}
where in the last inequality we have used  the fact that $\partial_t^2v=-\int_0^y\partial_t^2\partial_x\cdot u(t,x,\tilde y)d\tilde y$  (see Appendix \ref{secap} for details).
By using the argument  in  Section \ref{secderivation}, we can derive   estimates on $U^\eps$ and $\eps V^\eps$   as those given in \eqref{+pre}. 
Hence, by using the estimates in Theorems \ref{thhyp} and \ref{thmhas},  we   conclude that 
\begin{eqnarray*}
\begin{aligned}
	&\sup_{t \geq 0} \big(\abs{  U^\eps(t)}_{X_{\rho(t)/2}}^2 +\abs{\eps V^\eps(t)}_{X_{\rho(t)/2}}^2 \big)+\int_0^{+\infty}\big(\abs{  U^\eps(s)}_{X_{\rho(s)/2}}^2 +\abs{\eps V^\eps(s)}_{X_{\rho(s)/2}}^2 \big)ds\\
	&\leq  C \big(\abs{u_0^\eps-u_0}_{X_{\rho_0}}^2+\abs{u_1^\eps-u_1}_{X_{\rho_0}}^2 +\eps\abs{ v_0^\eps-v_0 }_{X_{\rho_0}}^2 +\eps\abs{ v_1^\eps-v_1 }_{X_{\rho_0}}^2 \big) \\
	&\quad +\eps^2C \int_0^{+\infty}\inner{ \abs{u}_{X_\rho} + \abs{u}_{X_\rho}^2+ \abs{\partial_tu}_{X_{\rho}}}^2ds \\
	&  \leq  C\inner{\abs{u_0^\eps-u_0}_{X_{2\rho_0}}^2+\abs{u_1^\eps-u_1}_{X_{2\rho_0}}^2}+ \eps^2C,
	\end{aligned}
\end{eqnarray*}
where we have  used  Theorem  \ref{thhyp} in the last inequality.  
Then this completes the proof of  Theorem \ref{thmlim}.

\bigskip
\noindent {\bf Acknowledgements.}
 The research of Wei-Xi Li  was supported by NSFC (Nos. 11961160716, 11871054, 12131017) and  the Natural Science Foundation of Hubei Province (2019CFA007).   And the research of Tong Yang
 was supported by the General Research Fund of Hong Kong CityU No. 11302020.

\appendix

\section{Some computation}\label{secap}

We now present the proof of  the estimates \eqref{st1}-\eqref{vsot} and \eqref{eforpressure}-\eqref{trau}.

\begin{proof}[{\bf Proof of \eqref{st1} and \eqref{st2}}]  In the proof $C$  denotes a 
generic constant  depending only on $\rho_0$ and the Sobolev embedding constant but independent of $\eps$.   To estimate the first term on the right of   
	\begin{eqnarray*}
\begin{aligned}
	R_{\eps} = -(U^\eps\cdot\partial_x)u^\eps- (u\cdot\partial_x)U^\eps  -V^\eps\partial_yu^\eps-v\partial_yU^\eps,
	\end{aligned}
\end{eqnarray*} 
we use the estimate 
\begin{eqnarray*}
\forall\ k\geq 0, \quad 	\frac{m!}{k!(m-k)!}\frac{L_{\rho/2,m}}{ L_{\rho,k+3}L_{\rho/2,m-k}}  \leq 2^{7-k} (k+3)^6  \rho^{-4}, 
\end{eqnarray*}
to compute 
\begin{equation}\label{appco}
\begin{aligned}
&	\sum_{m=0}^{+\infty}L_{\rho/2,m}^2	\big|\big(\partial_{x_j}^m \big((U^\eps\cdot\partial_x)u^\eps\big),\   \partial_{x_j}^m U^\eps+2\partial_t\partial_{x_j}^m U^\eps \big)_{L^2}\big|\\
& \leq    \sum_{m=0}^{+\infty} \sum_{k=0}^{m}   \frac{m!}{k!(m-k)!}\frac{L_{\rho/2,m}}{ L_{\rho,k+3}L_{\rho/2,m-k}}      L_{\rho,k+3
}\norm{ \partial_x \partial_{x_j}^k u^\eps}_{L^\infty} \\
	 &\qquad\qquad \qquad\quad\times  \Big(L_{\rho/2,m-k}\norm{   \partial_{x_j}^{m-k}U^\eps}_{L^2}\Big)\times \Big( L_{\rho/2,m} \norm{\partial_{x_j}^m U^\eps+2\partial_t\partial_{x_j}^m U^\eps}_{L^2}\Big)\\
& \leq  C\rho^{-4}\ \bigg[  \sum_{m=0}^{+\infty} \bigg(\sum_{k=0}^{m}  2^{-k} k^6    L_{\rho,k+2}\norm{ \partial_x \partial_{x_j}^k u^\eps}_{L^\infty} L_{\rho/2,m-k}\norm{   \partial_{x_j}^{m-k}U^\eps}_{L^2}\bigg)^2\,\bigg]^{1/ 2} \abs{  U^\eps}_{X_{\rho/2}}\\
	&\leq \frac{ C}{\rho^4}   \sum_{k=0}^{\infty}    2^{-k}k^6      L_{\rho,k+3}\norm{  \partial_x\partial_{x_j}^k  u^\eps}_{L^\infty} \bigg[\sum_{m=0}^{+\infty}       L_{\rho/2, m}^2\norm{    \partial_{x_j}^{ m}U^\eps}_{L^2}^2  \bigg]^{1/2}  \abs{U^\eps}_{X_{\rho/2}}\\
	&\leq \frac{ C}{\rho^4}  \abs{u^\eps}_{X_{\rho}} \abs{U^\eps}_{X_{\rho/2}}^2\leq  C  \abs{u^\eps}_{X_{\rho}} \abs{U^\eps}_{Y_{\rho/2}}^2,
	\end{aligned}
\end{equation}
where  we have used Young's inequality for   discrete convolution in the third inequality and the fact that $ \rho_0/2\leq \rho\leq \rho_0$ in the last line.  Similarly, 
\begin{equation*}
\begin{aligned}
 	\sum_{m=0}^{+\infty}L_{\rho/2,m}^2	\big|\big(\partial_{x_j}^m \big(v\partial_yU^\eps\big),\   \partial_{x_j}^m U^\eps+2\partial_t\partial_{x_j}^m U^\eps \big)_{L^2}\big| \leq  C  \abs{u}_{X_{\rho}} \abs{U^\eps}_{Y_{\rho/2}}^2.
	\end{aligned}
\end{equation*}
As for the third term in $R_\eps$,  we use  the fact  $V^\eps=-\int_0^y\partial_x\cdot U^\eps (t,x,\tilde y)d\tilde y $  to write
\begin{equation*}
	\sum_{m=0}^{+\infty}L_{\rho/2,m}^2	\big|\big(\partial_{x_j}^m (V^\eps\partial_yu^\eps),\   \partial_{x_j}^m U^\eps+2\partial_t\partial_{x_j}^m U^\eps \big)_{L^2}\big|\leq J_1+J_2
\end{equation*}
with
\begin{eqnarray*}
	\begin{aligned}
	J_1&=  \sum_{m=0}^{+\infty} \sum_{k=0}^{[m/2]}   \frac{m!}{k!(m-k)!}\frac{L_{\rho/2,m}}{ L_{\rho/2,k+2}L_{\rho/2,m-k+1}}      \Big(L_{\rho/2,k+2}\norm{  \partial_{x_j}^k \partial_yu^\eps}_{L_x^\infty L^2_y}\Big)\\
	 &\qquad\qquad \qquad\quad\times  \Big(L_{\rho/2,m-k+1}\norm{  \partial_x \partial_{x_j}^{m-k}U^\eps}_{L^2}\Big)\times \Big( L_{\rho/2,m} \norm{\partial_{x_j}^m U^\eps+2\partial_t\partial_{x_j}^m U^\eps}_{L^2}\Big),\\
	J_2 &= \sum_{m=0}^{+\infty} \sum_{k=[m/2]+1}^m   \frac{m!}{k!(m-k)!}\frac{L_{\rho/2,m}}{ L_{\rho,k+2}L_{\rho/2,m-k+1}}      \Big(L_{\rho,k+2}\norm{  \partial_{x_j}^k \partial_yu^\eps}_{L_x^\infty L^2_y}\Big)\\
	 &\qquad\qquad \qquad\quad\times  \Big(L_{\rho/2,m-k+1}\norm{ \partial_x  \partial_{x_j}^{m-k}U^\eps}_{L^2}\Big)\times \Big( L_{\rho/2,m} \norm{\partial_{x_j}^m U^\eps+2\partial_t\partial_{x_j}^m U^\eps}_{L^2}\Big).
	 \end{aligned}
\end{eqnarray*}
Repeating the computation in \eqref{fre} yields
\begin{eqnarray*}
	J_1 \leq   C  \abs{u^\eps}_{X_{\rho/2}}\abs{U^\eps }_{Y_{\rho/2}}^2.
\end{eqnarray*}
Moreover, by  observing
\begin{eqnarray*}
	\forall\  k\geq[m/2]+1,\quad     \frac{m!}{k!(m-k)!}\frac{L_{\rho/2,m}}{ L_{\rho,k+2}L_{\rho/2,m-k+1}} \leq 2^{7-k} (k+2)^6  \rho^{-4},
\end{eqnarray*}
and by  using  a similar argument as in  \eqref{appco},  we conclude
\begin{eqnarray*}
	J_2\leq  C    \abs{u^\eps}_{X_{\rho}}\abs{U^\eps }_{Y_{\rho/2}}^2.
\end{eqnarray*}
Thus,
\begin{eqnarray*}
	\sum_{m=0}^{+\infty}L_{\rho/2,m}^2	\big|\big(\partial_{x_j}^m (V^\eps\partial_yu^\eps),\   \partial_{x_j}^m U^\eps+2\partial_t\partial_{x_j}^m U^\eps \big)_{L^2}\big|\leq C    \abs{u^\eps}_{X_{\rho}}\abs{U^\eps }_{Y_{\rho/2}}^2.
\end{eqnarray*}
Similarly, we have
\begin{eqnarray*}
	\sum_{m=0}^{+\infty}L_{\rho/2,m}^2	\big|\big(\partial_{x_j}^m ((u\cdot\partial_x)U^\eps),\   \partial_{x_j}^m U^\eps+2\partial_t\partial_{x_j}^m U^\eps \big)_{L^2}\big|\leq C\abs{u}_{X_{\rho}}\abs{U^\eps }_{Y_{\rho/2}}^2.
\end{eqnarray*}
In summary, we have \eqref{st1}. 

Next we prove \eqref{st2}. Recall 
\begin{eqnarray*}
	S_\eps  =-\eps^2(U^\eps\cdot\partial_x)v^\eps-\eps^2(u\cdot\partial_x)V^\eps-\eps^2V^\eps\partial_yv^\eps-\eps^2v\partial_yV^\eps.
\end{eqnarray*}
Following the above argument used in \eqref{appco} with slight modification,  we can show that 
\begin{multline*}
	\sum_{m=0}^{+\infty}L_{\rho/2,m}^2	\big|\big(\partial_{x_j}^m \big(\eps^2(U^\eps\cdot\partial_x)v^\eps\big),\   \partial_{x_j}^m V^\eps+2\partial_t\partial_{x_j}^m V^\eps \big)_{L^2} \big|\\
	\leq C  \abs{\eps v^\eps}_{X_{\rho}} \abs{U^\eps}_{X_{\rho/2}} \abs{\eps V^\eps}_{X_{\rho/2}}  \leq  C  \abs{\eps v^\eps}_{X_{\rho}} \abs{U^\eps}_{Y_{\rho/2}}^2+C \abs{\eps v^\eps}_{X_{\rho}} \abs{\eps V^\eps}_{Y_{\rho/2}}^2.
\end{multline*} 
Similarly, by observing $\partial_yv^\eps=-\partial_x\cdot u^\eps,$ we have
\begin{eqnarray*}
	\sum_{m=0}^{+\infty}L_{\rho/2,m}^2	\big|\big(\partial_{x_j}^m \big(\eps^2V^\eps\partial_yv^\eps\big),\   \partial_{x_j}^m V^\eps+2\partial_t\partial_{x_j}^m V^\eps \big)_{L^2} \big| 
	\leq C  \abs{u^\eps}_{X_{\rho}}  \abs{\eps V^\eps}_{X_{\rho/2}}^2  \leq  C  \abs{u^\eps}_{X_{\rho}} \abs{\eps V^\eps}_{Y_{\rho/2}}^2. 
\end{eqnarray*}
Using  a similar computation as in \eqref{fre} and \eqref{appco} gives
\begin{eqnarray*}
		\sum_{m=0}^{+\infty}L_{\rho/2,m}^2	\big|\big(\partial_{x_j}^m \big(\eps^2(u\cdot\partial_x)V^\eps\big),\   \partial_{x_j}^m V^\eps+2\partial_t\partial_{x_j}^m V^\eps \big)_{L^2} \big| \leq   C \abs{u}_{X_{\rho}} \abs{\eps V^\eps}_{Y_{\rho/2}}^2.
\end{eqnarray*}
Finally,
\begin{eqnarray*}
		\sum_{m=0}^{+\infty}L_{\rho/2,m}^2	\big|\big(\partial_{x_j}^m \big(\eps^2v\partial_yV^\eps\big),\   \partial_{x_j}^m V^\eps+2\partial_t\partial_{x_j}^m V^\eps \big)_{L^2} \big| \leq   C \abs{u}_{X_{\rho}} \abs{\eps V^\eps}_{Y_{\rho/2}}^2.
\end{eqnarray*}
 Combining the above estimates yields \eqref{st2}. 
\end{proof}

\begin{proof}
	[{\bf Proof of \eqref{vsot}}] 
Observe that $\partial_t^2v=-\int_0^y\partial_t^2 \partial_x\cdot ud\tilde y$. Then it follows from  Definition \ref{defgev} that
	\begin{multline*}
		\sum_{m=0}^{+\infty}L_{\rho/2,m}^2\norm{ \partial_{x_j}^m  \partial_t^2 v}_{L^2}^2\leq \sum_{m=0}^{+\infty} \frac{L_{\rho/2,m}^2 }{L_{\rho,m+1}^2 } L_{\rho,m+1}^2  \norm{ \partial_{x_j}^m  \partial_t^2 v}_{L^2}^2\\
		\leq C \sum_{m=0}^{+\infty} \frac{m+1}{2^{ m+1}\rho}  L_{\rho,m+1}^2  \norm{ \partial_x\partial_{x_j}^m  \partial_t^2 u}_{L^2}^2\leq C\abs{\partial_t u}_{X_\rho}^2.
	\end{multline*}
	Similarly, using again the fact that $\partial_yv=-\partial_x\cdot u$, we have
	\begin{multline*}
		\sum_{m=0}^{+\infty}L_{\rho/2,m}^2\inner{\norm{ \partial_{x_j}^m  \partial_t v}_{L^2}^2+\norm{ \partial_{x_j}^m  \Delta_x v}_{L^2}^2+\norm{ \partial_{x_j}^m  \partial_y^2 v}_{L^2}^2} \\
		\leq  C  \sum_{m=0}^{+\infty}    L_{\rho,m+1}^2 \inner{\norm{ \partial_x\partial_{x_j}^m  \partial_t u}_{L^2}^2 +\norm{\partial_x \partial_{x_j}^m\partial_y u}_{L^2}^2}+  C\sum_{m=0}^{+\infty}    L_{\rho,m+3}^2  \norm{\partial_x \partial_{x_j}^m  \Delta_x u}_{L^2}^2  \leq C\abs{u}_{X_\rho}^2.
	\end{multline*}
	Combining the above estimates gives
	\begin{multline*}
        \sum_{m=0}^{+\infty} L_{\rho/2,m}^2	\big(-\eps^2\partial_{x_j}^m \big(\partial_t^2+\partial_t    -\eps^2\Delta_x-\partial_y^2\big) v,\   \partial_{x_j}^m V^\eps+2\partial_t\partial_{x_j}^m V^\eps \big)_{L^2} \\  \leq C \eps \inner{\abs{u}_{X_\rho}+\abs{\partial_t u}_{X_\rho}}\abs{\eps V^\eps}_{X_{\rho/2}}.  
	\end{multline*}
Moreover,  using a similar computation as in \eqref{appco} yields 
	\begin{equation*}
			\sum_{m=0}^{+\infty}L_{\rho/2,m}^2	\big(-\eps^2\partial_{x_j}^m \big( u\cdot\partial_x+v\partial_y  \big) v,\   \partial_{x_j}^m V^\eps+2\partial_t\partial_{x_j}^m V^\eps \big)_{L^2}\leq  C \eps  \abs{u}_{X_\rho}^2\abs{\eps V^\eps}_{X_{\rho/2}}.
	\end{equation*}
	Thus the estimate \eqref{vsot} follows. 
	\end{proof}
  
  \begin{proof}
  	[{\bf Proof of \eqref{eforpressure} and \eqref{trau}}]  Letting $t=0$ in \eqref{ellpre} and then applying $\partial_{x_j}^m$ to the both sides of the  equation, we obtain by standard elliptic theory that 
  	\begin{eqnarray*}
  		\begin{aligned}
  			\norm{\partial_{x_j}^m\partial_x p|_{t=0}}_{L_x^2}^2 \leq \inner{ \norm{\partial_{x_j}^m \big[(u_0\cdot \partial_x) u_0 +  ( \partial_x\cdot u_0) u_0 \big]}_{L^2}+  \norm{\partial_{x_j}^m\partial_y^2  u_0}_{L^2}}\norm{\partial_{x_j}^m\partial_x p|_{t=0}}_{L_x^2}.
  		\end{aligned}
  	\end{eqnarray*}
  	Then using a similar computation as in \eqref{appco} gives 
  	\begin{eqnarray*}
  		\begin{aligned}
  			\sum_{m=0}^\infty L_{\rho_0, m}^2  \norm{\partial_{x_j}^m\partial_x p|_{t=0}}_{L_x^2}^2\leq C \inner{ \abs{u_0}_{X_{2\rho_0}}^2+\abs{\partial_yu_0}_{X_{\rho_0}} }\Big(\sum_{m=0}^\infty L_{\rho_0, m}^2  \norm{\partial_{x_j}^m\partial_x p|_{t=0}}_{L_x^2}^2\Big)^{1/2}.
  		\end{aligned}
  	\end{eqnarray*}
  	Thus \eqref{eforpressure} follows. Similar argument holds for \eqref{trau}.
  \end{proof}


\end{document}